\documentclass[reqno,10pt]{amsart}
\usepackage{hyperref}

\usepackage{amsmath}
\usepackage{amsthm}
\usepackage{amssymb}
\usepackage{mathrsfs}
\usepackage{comment}
\usepackage{dsfont}
\usepackage{amssymb}
\usepackage{amsfonts}
\usepackage{mathpazo}
\usepackage{color}
\usepackage{xcolor}
\usepackage{yfonts}
\usepackage{paralist}
\usepackage{amsxtra}

\usepackage[normalem]{ulem}

\newtheorem{theorem}{Theorem}[section]
\newtheorem{corollary}[theorem]{Corollary}
\newtheorem{lemma}[theorem]{Lemma}
\newtheorem{conjecture}[theorem]{Conjecture}
\newtheorem{proposition}[theorem]{Proposition}
\newtheorem{prop}[theorem]{Proposition}
\newtheorem{cor}[theorem]{Corollary}
\theoremstyle{definition}

\newtheorem{claim}{Claim}
\newtheorem{definition}[theorem]{Definition}
\newtheorem{remark}[theorem]{Remark}

\numberwithin{equation}{section}

\newcommand{\NN}{\mathbb{N}}
\newcommand{\ZZ}{\mathbb{Z}}
\newcommand{\RR}{\mathbb{R}}
\newcommand{\CC}{\mathbb{C}}
\newcommand{\HH}{\mathbb{H}}

\newcommand{\cL}{\mathcal{L}}

\newcommand{\st}{\mathrm{st}}
\newcommand{\Diff}{\mathrm{Diff}}

\newcommand{\Stab}{\mathrm{Stab}}
\newcommand{\spann}{\mathrm{span}}

\newcommand{\Lie}{\mathrm{Lie}}
\newcommand{\Ad}{\mathrm{Ad}}
\newcommand{\ad}{\mathrm{ad}}
\newcommand{\htt}{\mathrm{ht}}

\newcommand{\SL}{\mathrm{SL}}
\newcommand{\SO}{\mathrm{SO}}
\newcommand{\OO}{\mathrm{O}}
\newcommand{\SU}{\mathrm{SU}}
\newcommand{\Sp}{\mathrm{Sp}}
\newcommand{\Lg}{\mathfrak{g}}
\newcommand{\Lh}{\mathfrak{h}}
\newcommand{\Lk}{\mathfrak{k}}
\newcommand{\Lp}{\mathfrak{p}}
\newcommand{\Lq}{\mathfrak{q}}
\newcommand{\La}{\mathfrak{a}}
\newcommand{\Lb}{\mathfrak{b}}

\newcommand{\Lc}{\mathfrak{c}}
\newcommand{\Ln}{\mathfrak{n}}
\newcommand{\Lu}{\mathfrak{u}}

\newcommand{\Lgl}{\mathfrak{gl}}
\newcommand{\Lsl}{\mathfrak{sl}}
\newcommand{\Lso}{\mathfrak{so}}
\newcommand{\Lsu}{\mathfrak{su}}
\newcommand{\Lsp}{\mathfrak{sp}}
\newcommand{\ra}{\rangle}
\newcommand{\sm}{\smallsetminus}
\newcommand{\mt}{\mapsto}

\def\eps{\varepsilon}

\def\reg{     1+\kappa}

\def\a{          \alpha}

\def\cF{          \mathcal F}
\def\cG{          \mathcal G}
\def\cH{          \mathcal H}
\def\cN{          \mathcal N}

\def\cW{          \mathcal W}
\def\scrV{          \mathscr V}

\def\supp{   {\rm supp}}

\def \R{{\mathbb R}}

\def \Z{{\mathbb Z}}
\def \N{{\mathbb N}}

\def \a{{\alpha}}
\def \ta{{\tilde\alpha}}

\newcommand{\prf}{{\begin{proof}}}
	\newcommand{\epf}{{\end{proof}}}

\newcommand{\ba}{{\mathbf a}}

\newcommand{\bu}{{\mathbf u}}
\newcommand{\bv}{{\mathbf v}}

\newcommand{\bh}{\mathbf{h}}
\newcommand{\bx}{\mathbf{x}}
\newcommand{\by}{\mathbf{y}}

\newcommand{\cE}{{\mathcal E}}

\newcommand{\calE}{{\mathcal E}}

\newcommand{\liealgg}{{\mathfrak g}}
\newcommand{\liealga}{{\mathfrak a}}

\renewcommand{\le}{\leqslant}
\renewcommand{\leq}{\leqslant}
\renewcommand{\ge}{\geqslant}
\renewcommand{\geq}{\geqslant}

\newcommand{\ary}{\begin{eqnarray}}
\newcommand{\eary}{\end{eqnarray}}
\newcommand{\aryst}{\begin{eqnarray*}}
\newcommand{\earyst}{\end{eqnarray*}}
\newcommand{\enmt}{\begin{enumerate}}
\newcommand{\eenmt}{\end{enumerate}}

\begin{document}
	
\title{Zimmer's conjecture for non-split semisimple Lie groups}

\author{Jinpeng An}
\address{Jinpeng An, School of Mathematical Sciences, Peking University, Beijing, 100871, China}
\email{anjinpeng@gmail.com}

\author{Aaron Brown}
\address{Aaron Brown, Northwestern University, Evanston, IL 60208, USA}
\email{awb@northwestern.edu}

\author{Zhiyuan Zhang}
\address{Zhiyuan Zhang, Department of Mathematics, Imperial College London, London SW7 2AZ, United Kingdom}
\email{zhiyuan.zhang@imperial.ac.uk}

\begin{abstract}
We prove many new cases of Zimmer's conjecture for actions by lattices in non-$\R$-split semisimple Lie groups $G$.  By prior arguments,  Zimmer's conjecture reduces to studying certain probability measures invariant under a minimal parabolic subgroup for the induced $G$-action.

Two techniques are introduced to give lower bounds on the dimension of a manifold $M$ admitting a non-isometric action.  First, when the Levi component of the stabilizer of the measure has higher-rank simple factors, cocycle superrigidity provides a lower bound on the dimension of $M$.  Second, when certain fiberwise coarse Lyapunov distributions are one-dimensional, a measure rigidity argument provides additional invariance of the measure if the associated root spaces are higher-dimensional.
\end{abstract}

\maketitle

%
	\tableofcontents


\section{Introduction}\label{sec Introduction}

As is well known, lattices in higher rank semisimple Lie groups exhibit various kinds of rigidity. This can be interpreted as the phenomenon that a lattice $\Gamma$ carries information about the ambient Lie group $G$. For instance, Margulis' superrigidity theorem states, roughly speaking, that any representation of $\Gamma$ can be extended to a representation of $G$ up to a bounded error. In particular, if we let $n(G)$ denote the minimal dimension of nontrivial real representations of the Lie algebra of $G$, then the dimension of any real representation $\rho$ of $\Gamma$ with unbounded image is at least $n(G)$.

Note that $\rho$ has unbounded image if and only if it does not preserve any inner product on the representation space, which can be viewed as a ``nontriviality condition" on $\rho$. If we change the interpretation of ``nontriviality" by only requiring $\rho(\Gamma)$ to be infinite, Margulis' superrigidity theorem with compact codomains (see, e.g., \cite[Section 2.3]{BFH}) further implies that the dimension of a real representation of $\Gamma$ with infinite image is at least $\min\{n(G),n(G_{\mathrm{cpt}})\}$, where $G_{\mathrm{cpt}}$ is the compact real form of  the complexification of the adjoint group $\Ad(G)$ of $G$.

Zimmer's conjecture is a collection of statements that extend the above rigidity results for linear representations to actions (namely, nonlinear representations) on compact manifolds. For a connected Lie group $G$, let $v(G)$ denote the minimal dimension of compact manifolds that admit nontrivial $G$-actions. If $G$ is noncompact semisimple, such manifolds are always covers of $G/Q$ where $Q$ is a maximal proper parabolic subgroup (see \cite[Theorem 6.6]{St}). Zimmer's conjecture predicts that if $G$ is semisimple of higher rank and $\Gamma\subset G$ is a lattice, then the minimal dimension of compact manifolds $M$ that admit $\Gamma$-actions without invariant Riemannian metrics coincides with $v(G)$.
In other words, if we regard nonlinear representations with invariant Riemannian metrics as ``trivial", then the minimal dimension of ``nontrivial" nonlinear representations of $\Gamma$ carries information about the same quantity for $G$.
This is a nonlinear extension of the above special case of Margulis' superrigidity theorem.
Strengthening the interpretation of ``triviality" of a $\Gamma$-action, Zimmer's conjecture also predicts that the minimal dimension of compact manifolds $M$ which admit actions $\alpha:\Gamma\to\Diff(M)$ with infinite images $\alpha(\Gamma)$ is at least $\min\{v(G),v(G_{\mathrm{cpt}})\}$. Let us rephrase the two conjectural statements as follows.

\begin{conjecture}[Zimmer's Conjecture, non-volume-preserving case]\label{C:Zimmer1}
Let $G$ be a connected semisimple Lie group with finite center, all of whose simple factors have real rank at least $2$. Let $\Gamma\subset G$ be a lattice, $M$ be a compact manifold, and $\alpha:\Gamma\to\Diff(M)$ be an action. Then the following statements hold:
\begin{itemize}
  \item[(1)] If $\dim(M)<v(G)$, then $\alpha$ preserves a Riemannian metric on $M$.
  \item[(2)] If $\dim(M)<\min\{v(G),v(G_{\mathrm{cpt}})\}$, then the image $\alpha(\Gamma)$ is finite.
\end{itemize}
\end{conjecture}

\begin{remark}
\begin{itemize}
  \item[(1)]   As noted above, the quantity $v(G)$ coincides with the minimal codimension of proper parabolic subgroups of $G$. On the other hand,
  $v(G_{\mathrm{cpt}})$ coincides with the minimal codimension of proper closed subgroups of $G_{\mathrm{cpt}}$ (which is denoted by $d(G)$ in \cite{BFH,BFH-new}). Both quantities depend only on the Lie algebra of $G$.
  \item[(2)] In view of Margulis' superrigidity theorem with compact codomains, statement (1) in Conjecture \ref{C:Zimmer1} implies statement (2) (see \cite[Section 2.3]{BFH}).
\end{itemize}
\end{remark}

For the case that $G$ is $\RR$-split, Conjecture \ref{C:Zimmer1} has been proved in the recent breakthrough \cite{BFH,BFH2,BFH-new}.
For non-split groups, a weaker result has also been proved in \cite{BFH,BFH-new}, namely, Conjecture \ref{C:Zimmer1} holds if the quantity $v(G)$ is replaced by the so-called minimal resonant codimension $r(G)$. For example, if $G$ is complex semisimple, one has $r(G)=v(G)/2$. The quantities $v(G)$, $v(G_{\mathrm{cpt}})$ and $r(G)$ for some simple Lie groups $G$ are listed in Appendix \ref{S:vnd}. See also \cite{Bro,BDZ,Can,Fis1,Fis2,Fis3} for surveys and related works.

This paper is concerned with the non-split case. One of the main results is as follows.

\begin{theorem}\label{T:main-1}
Let $G$ be one of the following Lie groups:
\begin{itemize}
  \item[(i)] A connected complex semisimple Lie group without simple factors of rank $1$;
  \item[(ii)] $\SL(n,\HH)$ with $n\ge 5$;
  \item[(iii)] $\SO^+(m,n)$ with either $(m,n)=(4,2)$ or $3\le n<m\le \frac12(n^2+n+4)$;
  \item[(iv)] $\SU(m,n)$ with $2\le n \le m\le \frac12(n^2-n+2)$;
  \item[(v)] $\SO^*(2n)$ with $n\ge 14$;
  \item[(vi)]  $E II$ (following the notation in \cite[App.\ C.4]{Kn}).
\end{itemize}
Let $\Gamma\subset G$ be a lattice, $M$ be a compact manifold, and $\alpha:\Gamma\to\Diff^{\reg}(M)$ be an action. Then the following statements hold:
\begin{itemize}
  \item[(1)] If $\dim(M)<v(G)$, then $\alpha$ preserves a continuous Riemannian metric.
  \item[(2)] If $\dim(M)<\min\{v(G),v(G_{\mathrm{cpt}})\}$, then $\alpha(\Gamma)$ is finite.
\end{itemize}
\end{theorem}

For $G=\Sp(m,n)$ with $6\le n\le m\le \frac14(n^2-3n+6)$, we will prove a weaker result where $v(G)$ is replaced by $v(G)-2$ (see Theorem \ref{T:Sp}). 
For $G=\SL(n,\CC)$ and $\Gamma$ is cocompact, Theorem \ref{T:main-1} was proved by the third-named author in the preprint \cite{Zh} (which will not be submitted for publication). 

Let us also state the volume-preserving analogue of Conjecture \ref{C:Zimmer1}.

\begin{conjecture}[Zimmer's Conjecture, volume-preserving case]\label{C:Zimmer2}
Let $G, \Gamma, M$ and $\alpha$ be as in Conjecture \ref{C:Zimmer1}, and suppose that the action $\alpha$ preserves a volume density on $M$. Then the following statements hold:
\begin{itemize}
  \item[(1)] If $\dim(M)<n(G)$, then $\alpha$ preserves a Riemannian metric on $M$.
  \item[(2)] If $\dim(M)<\min\{n(G),v(G_{\mathrm{cpt}})\}$, then the image $\alpha(\Gamma)$ is finite.
\end{itemize}
\end{conjecture}

\begin{remark}
\begin{itemize}
  \item[(1)] Recall that $n(G)$ is the minimal dimension of nontrivial real representations of the Lie algebra of $G$.
One always has $n(G)>v(G)$ (see \cite[Lemma 6.2]{St}).
  \item[(2)] As in the case of Conjecture \ref{C:Zimmer1}, item (1) in Conjecture \ref{C:Zimmer2} implies item (2).
\end{itemize}
\end{remark}

For volume-preserving actions, we will prove:

\begin{theorem}\label{T:main-2}
Suppose that $G$ is (isogenous to) one of the following groups:
\begin{itemize}
  \item A connected complex semisimple Lie group without simple factors of rank $1$;
  \item $\SO^+(n,n+2)$  or $\SU(n,n)$ with $n\ge2$;
  \item  $E II$ (following the notation in \cite[App.\ C.4]{Kn}).
\end{itemize}
Let $\Gamma$ be a lattice in $G$, $M$ be a compact manifold with $\dim(M)=v(G)$, and $\alpha:\Gamma\to\Diff^{\reg}(M)$ be an action preserving a smooth volume density. Then $\alpha$ preserves a continuous Riemannian metric.
\end{theorem}

\subsection{Outline of new arguments}
The broad outline of proof for the results in this paper roughly follow those from \cite{BFH,BFH-new}.  Namely, we consider a higher-rank semisimple
Lie group $G$, a lattice $\Gamma\subset G$, and an action $\alpha\colon \Gamma\to \Diff^{r}(M)$ for $r>1$.  If the action $\alpha(\Gamma)$ is not isometric for any $C^0$ Riemannian metric then, for the induced action on the suspension space $M^\alpha$ (see Section \ref{sec:susspace} below), there is an ergodic, $A$-invariant probability measure $\mu$, projecting to the normalized Haar measure on $G/\Gamma$, with a non-zero fiberwise Lyapunov exponent.  A clever argument (Lemma \ref{L:PA} below) shows that we may always assume $\mu$ is further invariant under a minimal parabolic subgroup $P\subset G$.

Suppose that  $\mu$ were $G$-invariant.  Zimmer's cocycle superrigidity theorem (Theorem \ref{thm ZCS} below) applied to the fiberwise derivative cocycle over the $G$-action on $(M^\alpha, \mu)$ would then give a lower bound on dimension, $\dim (M)\ge n(G)$.
In general, $\mu$ will be $Q$-invariant for some (parabolic) $P\subset Q\subset G$.
Here we introduce two new techniques to provide lower bound on the dimension of $M$ by either (1) applying cocycle superrigidity to simple factors of Levi subgroup of $Q$ or (2) using measure rigidity arguments to produce extra invariance of $\mu$ assuming certain (coarse) fiberwise Lyapunov subspaces are 1-dimensional.  Both approaches provide lower bounds on $\dim (M)$ by providing lower bounds on the dimension of certain fiberwise Lyapunov subspaces.

\subsubsection{Cocycle superrigidity}
Let $Q$ be the stabilizer of the measure $\mu$ as above.  Let $L$ be a choice of (reductive) Levi subgroup in $Q$ containing a choice of maximal $\R$-split Cartan subgroup $A$ (relative to which we defined the restricted roots of the Lie algebra.)
When some simple factor of $L$ is of higher-rank, we may apply Zimmer's cocycle superrigidity theorem to fiberwise derivative cocycle over the action of $L$ on $(M^\alpha, \mu)$.  Lower bounds on the dimension of non-trivial real representations of the higher-rank simple factor and combinatorics determining $Q$ then give lower bounds on the dimension of the fiberwise Lyapunov subspaces for the $A$-action on $(M^\alpha, \mu)$, and hence lower bounds on the dimension of $M$.  See Proposition \ref{P:real-1}.

\subsubsection{Measure rigidity}
Again, let $Q$ be the stabilizer of the measure $\mu$ above. The Lie algebra $\mathfrak{q}$ of $Q$ is saturated by (coarse) root spaces $\mathfrak{g}_\beta$.  Let $\beta$ be a root of $G$  such that $\mathfrak{g}_\beta\cap \mathfrak{q} = \{0\}$ (we further assume $\frac 1 2 \beta$ is not a root).
It follows from \cite[Prop.\ 5.1]{BRHW} that there exists a fiberwise Lyapunov functional $\lambda^F$ that is positively proportional to $\beta$.  Let $\chi^F$ denote the associated coarse fiberwise Lyapunov functional and suppose the associated coarse Lyapunov subspace $E^{\chi^F}(x)$ is 1-dimensional for $\mu$-a.e.\ $x$.  (See Section \ref{sec:susspace} for definitions).  In this case, we prove a dichotomy: either $\dim (\Lg_{\beta}\oplus \Lg_{2\beta})= 1$ or $\mu$ is invariant under a 1-parameter subgroup of $G_{[\beta]}=\exp(\Lg_\beta\oplus\Lg_{2\beta} )$; since $\mu$ is assumed $P$-invariant it follows that $\mu$ is $G_{[\beta]}$-invariant.  This contradicts the definition of $Q$ and shows that if $\dim (\Lg_{\beta}\oplus \Lg_{2\beta})\ge 2$, the dimension of $E^{\chi^F}(x)$ is at least 2 for $\mu$-a.e.\ $x$.  See Theorem \ref{prop:dim2} below.

\section{Preliminaries} \label{sec preliminaries}

\subsection{Structure of semisimple Lie groups} \label{subsec Structure of semisimple Lie groups}

Throughout this paper, we let $G$ be a connected real semisimple Lie group with finite center and without simple factors of real rank less than $2$. Let $\Lg=\Lie(G)$ be the Lie algebra of $G$. We fix a Cartan decomposition $\Lg=\Lk\oplus\Lp$ and a maximal abelian subspace $\La$ of $\Lp$. Let $K$ and $A$ be the connected Lie subgroups of $G$ with Lie algebras $\Lk$ and $\La$, respectively. Then $K$ is a maximal compact subgroup, and $A$ is a maximal $\RR$-split Cartan subgroup of $G$. Let $C$ denote the centralizer of $A$ in $K$, and let $\Lc=\Lie(C)$.

By a slight abuse of the notations, for $\lambda \in \La^*$, $\ba\in\La$ and $a = \exp(\ba) \in A$, we will use $\lambda(a)$ to denote $\lambda(\ba)$. For $\lambda \in \La^*$, we denote $[\lambda] := \R_+ \lambda$ where $\RR_+=(0,\infty)$.

Let $\Sigma=\Sigma(\Lg,\La)\subset\La^*$ be the restricted root system of $\Lg$ relative to $\La$, and denote $\widehat{\Sigma}:=\{[\beta]:\beta\in\Sigma\}$.
Then we have the restricted root space decomposition $\Lg=\Lc\oplus\La\oplus\bigoplus_{\beta\in\Sigma}\Lg_\beta$. For every $[\beta]\in\widehat{\Sigma}$, the space $\Lg_{[\beta]}:=\bigoplus_{\beta'\in[\beta] \cap \Sigma}\Lg_{\beta'}$ is a subalgebra of $\Lg$.
We let $G_{[\beta]}$ denote the connected Lie subgroup of $G$ with Lie algebra $\Lg_{[\beta]}$.
Note that the only three possibilities of $[\beta] \cap \Sigma$ are $\{\beta\}$, $\{\beta, \beta/2\}$ and $\{\beta, 2\beta\}$. In the first two cases, namely when $2\beta \notin \Sigma$, $\Lg_\beta$ is also a subalgebra of $\Lg$, and we let $G_\beta$ denote the connected Lie subgroup of $G$ with Lie algebra $\Lg_\beta$. In the first case we have $G_{\beta} = G_{[\beta]}$, while in the second case $G_{\beta} \subsetneq G_{[\beta]}$.

We choose and fix a set of simple roots $\Pi\subset\Sigma$, and let $\Sigma^+$ denote the corresponding set of positive roots. For simplicity, we identify $\Pi$ with the set of vertices of the associated Dynkin diagram. Thus $\beta,\gamma\in\Pi$ are connected (by one or more edges) if and only if $\langle\beta,\gamma\ra<0$, where $\langle\cdot,\cdot\ra$ is the inner product on $\La^*$ induced from the Killing form of $\Lg$.

Let $\Ln=\bigoplus_{\beta\in\Sigma^+}\Lg_\beta$, and $N$ be the connected Lie subgroup of $G$ with Lie algebra $\Ln$. A closed subgroup $Q$ of $G$ is said to be parabolic if it contains a conjugate of the group $P:=CAN$. In particular, $Q$ is a minimal parabolic subgroup if and only if it is conjugate to $P$. We refer to parabolic subgroups containing $P$ as the standard ones.
If $Q\subset G$ is a standard parabolic subgroup, we denote
$$\Sigma_Q:=\{\beta\in\Sigma:\Lg_\beta\subset\Lie(Q)\}, \qquad \widehat{\Sigma}_Q:=\{[\beta]:\beta\in\Sigma_Q\}.$$
Then
$$\Lie(Q)=\Lc\oplus\La\oplus\bigoplus_{\beta\in\Sigma_Q}\Lg_\beta=\Lc\oplus\La\oplus\bigoplus_{[\beta]\in\widehat{\Sigma}_Q}\Lg_{[\beta]}.$$
There is a bijection $Q\mt\Pi_Q$ between standard parabolic subgroups $Q\subset G$ and subsets $\Pi_Q\subset\Pi$ such that $\Sigma_Q=\Sigma^+\cup(\spann(\Pi_Q)\cap\Sigma)$.

\subsection{Suspension space}\label{sec:susspace}

Let $\Gamma\subset G$ be a lattice.
Throughout this paper, we denote by $m_{G/\Gamma}$ the Haar measure on $G/\Gamma$.
Let $M$ be a compact manifold, and $\alpha:\Gamma\to \Diff^{\reg}(M)$ be a $C^{\reg}$ action.

As in the previous works \cite{BFH,BFH-new}, our proof of Theorems \ref{T:main-1} and \ref{T:main-2} will be based on investigation of the induced $G$-action on the suspension space $M^\alpha$. By definition, $M^\alpha$ is the space of $\Gamma$-orbits in the product space $G\times M$, where $\Gamma$ acts on $G\times M$ from the right by $(g,x).\gamma=(g\gamma,\alpha(\gamma)^{-1}x)$. The left $G$-action on $G\times M$ given by $g'.(g,x)=(g'g,x)$ commutes with the $\Gamma$-action, hence induces a $G$-action on $M^\alpha$, which will be denoted by $\tilde{\alpha}$.
The suspension space $M^\alpha$ is a fiber bundle over $G/\Gamma$ with fiber $M$, and the induced action $\tilde{\alpha}: G\to\Diff^{\reg}(M^\alpha)$ preserves the fiber foliation.

By an abuse of notation, for $\ba\in\La$, let us write $\tilde{\alpha}(\ba)$ for the expression $\tilde{\alpha}(\exp(\ba))$.
Let $\pi: M^\alpha\to G/\Gamma$ be the canonical projection,  and consider the fiberwise tangent bundle $E_F=\ker(d\pi)$ of $M^\alpha$. If $\mu$ is an $A$-invariant Borel probability measure on $M^\alpha$ such that for any $\ba \in \La$ and $1 \leq q < \infty$,
\ary  \label{eq integrability}
\int_{M^\alpha} \max(1, \log \| D_x \ta(\ba) \| )^q d\mu(x) < \infty,
\eary
we may define for any  $\ba\in\La$ the average top fiberwise Lyapunov exponent of $\tilde{\alpha}(\ba)$ with respect to $\mu$ to be
$$\lambda_+^F(\ba,\mu)=\liminf_{n\to\infty}\frac1n\int_{M^\alpha}\log\|D_x\tilde{\alpha}(n\ba)|_{E_F(x)}\|d\mu(x).$$

Recall that we have the following theorem (we state the version that we will apply in our setting).
\begin{theorem} \label{thm MET}
Assume that $\mu$ is an $A$-invariant and $A$-ergodic Borel probability measure on $M^\alpha$ satisfying \eqref{eq integrability},
then there exists:
\enmt
\item an integer $r \geq 1$ and a full measure, $A$-invariant subset $R_{reg} \subset M^{\a}$;
\item an $A$-invariant measurable family of linear functionals $\cL = \{ \lambda_i : \liealga \to \R \}_{1 \leq i \leq r}$;
\item and a family of mutually transverse, $D\tilde{\alpha}$-invariant, measurable subbundles $E_{\lambda_i} \subset TM^\a$ with $T_x M^\a = \bigoplus_{i=1}^r E_{\lambda_i}(x)$ for $x \in R_{reg}$
\eenmt
such that for all $v \in E_{\lambda_i}(x) \setminus \{0\}$
\ary \label{eq propertyregularpoint}
\lim_{\La\ni\ba\to\infty}\frac{\log\|D_x\tilde{\alpha}(\ba)(v)\|-\lambda_i(\ba)}{\|\ba\|}=0.
\eary
\end{theorem}
As usual, we will call $R_{reg}$ the set of {\it Lyapunov regular points}, call  each $\lambda_i$ a {\it Lyapunov functional}, and call $E_{\lambda_i}$ the corresponding {\it Oseledects subspace}.

Given $\epsilon > 0$ and $C > 0$, a {\it Pesin set } $R(\epsilon, C)$ is defined as the set of $x \in R_{reg}$ such that for any $\ba \in \liealga$, for $1 \leq i \leq r$, for $v \in E_{\lambda_i}(x)$, we have
\enmt
\item  $C^{-1} e^{\lambda_i(\ba) - \frac{1}{2}\epsilon \|\ba\|} \|v\| \leq \| D\tilde{\alpha}(\ba)(v)\| \leq C e^{\lambda_i(\ba) + \frac{1}{2}\epsilon \|\ba\|} \|v\|$;
\item $\angle(E_{\lambda_i}(x),  N(x)) \geq C^{-1}$, and $\angle(E_{\lambda_i}(x), E_{\lambda_j}(x) ) \geq C^{-1}$ for all $1 \leq i < j \leq r$.
\eenmt

Now fix an $A$-invariant and $A$-ergodic Borel probability measure $\mu$ on $M^\alpha$. We denote by $\cL(A, \mu)$ the set of linear functions in Theorem \ref{thm MET}.
Notice that $E_F$ is invariant under $D\tilde{\alpha}$, we deduce that $E_F$ is decomposed into a direct sum of Oseledects subspaces.
Then one define {\it fiberwise Lyapunov functionals}  $\cL_F(A,\mu)$ to be the set of Lyapunov functionals for the induced $A$-action on $E_F$ with respect to $\mu$. In other words, we have
$$E_F=\bigoplus_{\lambda\in\cL_F(A,\mu)}E_{F, \lambda}$$
where $E_{F, \lambda} = E_{\lambda} \cap E_F$.
We also make the convention that if $\lambda\in\La^*\sm\cL_F(A,\mu)$, then $E_{F, \lambda}$ is the zero bundle.
Note that for $\ba\in\La$ we have
$$\lambda_+^F(\ba,\mu)=\max_{\lambda\in\cL_F(A,\mu)}\lambda(\ba).$$

Let $P$ be the minimal standard parabolic subgroup of $G$.
If, instead, $\mu$ is a $P$-invariant $P$-ergodic Borel probability measure on $M^\alpha$ and $\{\mu^e_x:x\in M^\alpha\}$ is the $A$-ergodic components of $\mu$, then by \cite[Claim 5.2]{BRHW}, the set $\cL_F(A,\mu^e_x)$ is $\mu$-a.e. constant, which we denote again by $\cL_F(A,\mu)$. Moreover, for each $\lambda\in\cL_F(A,\mu)$, $\dim E_{F, \lambda}(\mu^e_x)$ is also $\mu$-a.e. constant.

We now introduce a collection of $\ta(\liealga)$-invariant foliations.
For $\ba \in \liealga$ and $x \in M^\a$, we define
\aryst
\cW^{-}_{\ta(\ba)}(x) = \{ y \in M^\a \mid \limsup_{n \to \infty} \frac{1}{n} \log d(\ta(n\ba)(x), \ta(n\ba)(y)) < 0 \}.
\earyst
From Pesin's theory we know that for $\mu$-a.e. $x\in M^\alpha$,  $\cW^{-}_{\tilde{\alpha}(\ba)}(x)$ is an embedded manifold whose tangent space at $x$ is $\bigoplus_{\lambda \in \cL_F(A,\mu): \lambda(\ba) < 0} E_{\lambda}(x)$.
We let $L$ be a subset of $\liealga \setminus \{0\}$ which freely generates a lattice in $\liealga$. Then for any $\lambda \in \cL(A, \mu)$, for $\mu$-a.e. $x$, we let $\cW^{[\lambda]}(x)$ be the path connected component at $x$ of
\aryst
\bigcap_{\ba \in L : \lambda(\ba) < 0} \cW^{-}_{\ta(\ba)}(x).
\earyst
For every $\lambda \in \cL(A, \mu)$,
for $\mu$-a.e. $x$, $\cW^{[\lambda]}(x)$ is an embedded manifold whose tangent space at $x$ is
\aryst
E^{[\lambda]}(x) = \bigoplus_{\lambda' \in \cL(A, \mu) \cap [ \lambda ] } E_{\lambda'}(x).
\earyst

We say that $\widehat\cL(A, \mu)$,  the $\R_+$ proportionality equivalent classes in $\cL(A, \mu)$,
is the collection of {\it coarse Lyapunov functionals}.
For each $\chi = [\lambda] \in \widehat\cL(A, \mu)$,
we say that $\cW^{\chi}(x)$ is the {\it coarse invariant manifold} at $x$, and $E^{\chi}
(x)$ is the {\it coarse invariant subspace} at $x$.  Similarly, we define
 {\it fiberwise coarse Lyapunov functionals}  $\widehat\cL_F(A, \mu)$  to be the $\R_+$ proportionality equivalent classes in $\cL_F(A, \mu)$, and denote
$\cW^{\chi}_F(x) = \cW^{\chi}(x) \cap \pi^{-1}(\pi(x))$ and $E^{\chi}_F(x) = E^{\chi}(x) \cap E_F$ for each $\chi \in \widehat\cL_F(A, \mu)$.
We let $\cW^{\chi}$, resp. $\cW^{\chi}_F$, denote the partition of $M^\a$ into coarse invariant manifolds  $\cW^{\chi}(x)$, resp. $\cW^{\chi}_F(x)$, $x \in M^\a$.
Note that the set of $G$-orbits form a foliation, denoted by $\cW_G$, that is transverse to every fiber of $\pi$. More generally, for each closed subgroup $H \subset G$, we denote by $\cW_H$ the orbital foliation of $H$. We will often consider $\cW_{G_{[\beta]}}$, which will be denoted by $\cW^{[\beta]}_{G}$, for some $\beta \in \Sigma$.

In applications, we need to make a slight generalization of the above construction.
Given a subgroup $A_0$ of $A$ with Lie algebra $\La_0 \subset \La$, we have Lyapunov functionals and their assocated Oseledects subspaces for the restriction $\ta|_{A_0}$.
In this case, a Lyapunov functional $\lambda$ belongs to $\La_0^*$, and its Oseledects subspace is
\aryst
E_{\lambda} = \bigoplus_{\substack{\lambda' \in  \cL(A, \mu): \\ \lambda'|_{\La_0} = \lambda }} E_{\lambda'}.
\earyst
We can define analogously  for $\ta|_{A_0}$ the coarse Lyapunov functionals $\widehat{\cL}(A_0, \mu)$ and the fiberwise coarse Lyapunov functionals $\widehat{\cL}_F(A_0, \mu)$.

In the rest of the paper, we will often consider some $A$-invariant measure on $M^{\alpha}$ which projects to the Haar measure on $G/\Gamma$. In this case we have the following.
\begin{lemma} \label{lem uniformexponential}
Fix some $x_0 \in M^{\a}$.
There exist constants $\kappa_0, C > 0$ depending only on $\alpha$ such that
 for any $A$-invariant, $A$-ergodic probability measure $\mu$ on $M^{\alpha}$ satisfying $\pi_*(\mu) = m_{G/\Gamma}$, we have $ \int e^{\tau_0 d(x_0, x)} d\mu(x)  < C$.
\end{lemma}
\begin{proof}
This is essentially contained in \cite{BFH-new}.
By the construction of the metric $d$ on $M^{\a}$, the fibers of $\pi$ are uniformly comparable. Then we have $d(x,y) - B \leq d_{G/\Gamma}(\pi(x), \pi(y)) \leq d(x,y)$.
As $\pi_*(\mu) = m_{G/\Gamma}$, the inequality in the lemma then follows from the fact that there exists $\tau_0 > 0$ such that
\aryst
\int e^{\tau_0 d(\pi(x_0), \pi(x))} d\mu(x) \leq \int e^{\tau_0 d_{G/\Gamma}(\pi(x_0), z)} dm_{G/\Gamma}(z) < \infty.
\earyst
The last inequality above is well-known.
\end{proof}

\subsection{Conditional measures along invariant foliations} \label{sec conditionalmeasure}

Let $\cF$ be an $A$-invariant foliation. For example, $\cF$ can be:
\enmt
\item $\cW^{\chi}$, $\cW^{\chi}_F$ for some $\chi \in \widehat\cL(A, \mu)$, or
\item The orbital foliation $\cW_{H}$ for some closed subgroup $H \subset G$, or
\item $\cW^{-}_{\ta(\ba)}$ for some $\ba \in \liealga$.
\eenmt
Given a family of measurable partition $\xi$ that is subordinate to $\cF$,
 we denote by $\sigma(\xi)$ the $\sigma$-algebra generated by $\xi$, and
there is a $\sigma(\xi)$-measurable assignment $x \mapsto \mu^{\xi}_x$, where $\mu^{\xi}_x$ is a probability measure on $\xi(x)$ for $\mu$-a.e. $x$, and for any measurable $X \subset M$ we have
\ary \label{eq measuredisintegration}
\mu(X) = \int \mu^{\xi}_x(X) d\mu(x).
\eary
We call $\mu^{\xi}_x$ the {\it conditional  measure} of $\mu$ at $x$ with respect to $\xi$.
Moreover, for another measurable partition $\xi'$ that is subordinate  to $\cF$, we have for $\mu$-a.e. $x$ that
\begin{align} \label{eq compatible measures}
\mu^{\xi}_x(\xi(x)\cap \xi'(x))^{-1} \mu^{\xi}_x |_{\xi(x)\cap \xi'(x)} = \mu^{\xi'}_x(\xi(x)\cap \xi'(x))^{-1}  \mu^{\xi'}_x |_{\xi(x)\cap \xi'(x)}.
\end{align}
This equality allows us to define a family of {\it leafwise measures} of $\mu$, denoted by $\{ [\mu^{\cF}_x]\}_{x \in M}$, where $\mu^{\cF}_x$ is a Radon measure on $\cF(x)$, defined up to a multiplicative constant, for $\mu$-a.e. $x$.
Since $\mu$ is $A$-invariant, we know that for every $\ba \in \liealga$, for $\mu$-a.e. $x$
\ary \label{eq ainvarianceconditionalmeasure}
[\ta(\ba)_*\mu^{\cF}_x] = [\mu^{\cF}_{\ta(\ba)(x)}].
\eary
Moreover, if $H$ is a closed subgroup of $G$, as $\ta$ is a locally free $G$-action, we can, for $\mu$-a.e. $x$, pullback $\mu^{\cW_H}_x$ to get a Radon measure $\mu^{H}_x$ on $H$,  which will be called the {\it conditional measure of $\mu$ along $H$ at $x$}.  We have the following proposition (see for example \cite[Proposition 5.7]{EK2}).
\begin{prop}\label{prop HinvariantHhaar}
The measure $\mu$ is $H$-invariant if and only if for $\mu$-a.e. $x$, the measure $\mu^H_x$ is a left-invariant Haar measure on $H$.
\end{prop}

We will also use the following result, which follows from the combination of Theorem 8.4 and Proposition 8.8 in \cite{EK2}.

\begin{prop}\label{lem noncommutinginvariance}
Let $\mu$ be an $A$-invariant Borel probability measure on $M^\alpha$, and let $\beta_1, \beta_2 \in \Sigma$ be linearly independent roots such that $\beta_1+\beta_2\in\Sigma$.
Then for $\mu$-a.e. $x\in M^\alpha$ and any $\bu_i \in \Lg_{[\beta_i]}$  with $\exp(\bu_i) \in \supp \mu^{G_{[\beta_i]}}_x$ \ $(i=1,2)$,
the conditional measure $\mu^{G_{[\beta_1+\beta_2]}}_x$ is invariant under the right translation of $\exp([\bv_1,\bv_2])$ on $G_{[\beta_1+\beta_2]}$, where $\bv_i\in\Lg_{\beta_i}$ is the projection of $\bu_i$ on $\Lg_{\beta_i}$ along $\bigoplus_{\beta\in[\beta_i] \cap \Sigma \sm\{\beta_i\}}\Lg_{\beta}$,
\end{prop}

\begin{proof}
Let $H$ be a ${\bf m}$-stable subgroup containing $G_{[\beta_1]}$ and $G_{[\beta_2]}$ for some ${\bf m} \in \La$ (see \cite[Page 174]{EK2} for definition).
We see that $\bv_1$ and $\bv_2$ are weight components of support vectors of $\mu$ as defined in \cite{EK2} (see  \cite[Page 221]{EK2}). Then we deduce from \cite[Prop 8.8]{EK2} that $\mu^{H}_x$ is invariant under the right translation of $\exp([\bv_1,\bv_2])$. We may conclude the proof by  \cite[Prop 8.4]{EK2} and by noting that $\exp([\bv_1,\bv_2]) \in G_{[\beta_1 + \beta_2]} \subset H$.
\end{proof}

\subsection{Entropy theory}

We recall some basic definitions and facts from the entropy theory (see for instance \cite{BRHW1}).

Let $\mu$ be a Borel probability measure on $M^\alpha$.
Given measurable partitions $\eta, \xi$ of $(M^{\alpha}, \mu)$, the mean conditional entropy of $\eta$ relative to $\xi$ is
\aryst
H_{\mu}(\eta \mid \xi) = \int - \log(\mu^{\xi}_x(\eta(x)))  d\mu(x).
\earyst
The entropy of $\eta$ is $H_{\mu}(\eta) = H_{\mu}(\eta \mid \{ \emptyset, M^{\alpha} \} )$. Note that if $H_{\mu}(\eta) < \infty$ then $\eta$ is necessarily countable.  Let $f: (M^{\alpha}, \mu) \to (M^{\alpha}, \mu)$ be an invertible, measurable, $\mu$-preserving transformation. Let $\eta$ be an arbitrary measurable partition of $(M^{\alpha}, \mu)$. We define
\aryst
\eta^{+} = \wedge_{i=0}^{\infty} f^i(\eta), \quad \eta^{f} = \wedge_{i \in \Z} f^i(\eta).
\earyst
We define the  entropy of $f$ conditioned on partition $\eta$ to be
\aryst
h_{\mu}(f, \eta) = H_{\mu}(\eta \mid f(\eta^+) ).
\earyst
We define the $\mu$-entropy of $f$ to be $h_{\mu}(f) = \sup\{  h_{\mu}(f, \eta) \}$ where the supremum is taken over all measurable partitions of $(M^{\alpha}, \mu)$.

Given a measurable partition $\eta$ of $(M^{\alpha}, \mu)$, we define the entropy of $f$ subordinate to $\eta$ to be
\aryst
h_{\mu}(f \mid \eta) = \sup \{ h_{\mu}(f, \xi) : \xi \geq \eta \}.
\earyst

Consider a $f$-invariant measurable foliation $\cF$ of $M^{\alpha}$. We say that a measurable partition $\eta$ of $(M^{\alpha}, \mu)$ is subordinate to $\cF$ if for a.e. $x \in M^{\alpha}$
\enmt
\item $\eta(x) \subset \cF(x)$,
\item $\eta(x)$ contains an open neighborhood of $x$ in $\cF(x)$,
\item $\eta(x)$ is precompact in $\cF(x)$ for a positive measure set of $x$.
\eenmt
Suppose in addition that $f$ is ergodic with respect to $\mu$.
We define the entropy of $f$ subordinate to $\cF$ is
\aryst
h_{\mu}(f \mid \cF) = \sup \{  h_{\mu}(f \mid \xi) :  \mbox{ $\xi$ is subordinate to $\cF$ } \}.
\earyst
More generally, if $\eta$ is any measurable partititon of $(M^{\alpha}, \mu)$ we define
\aryst
h_{\mu}(f \mid \eta \vee \cF) = \sup \{ h_{\mu}(f \mid \xi \vee \eta)  : \mbox{$\xi$ is subordinate to $\cF$}  \}.
\earyst

Using the above definitions, we can make sense of the notations such as $h_\mu( \ta( a ) )$,  $h_\mu(\ta( a ) \mid \cW^\chi)$  and $h_\mu( \ta( a ) \mid \cW^\chi_F)$  for an $A$-invariant $A$-ergodic invariant measure $\mu$ on $M^{\alpha}$, which we will abbreviate as $h_\mu(  a )$,  $h_\mu( a \mid \cW^\chi)$ and  $h_\mu(  a \mid \cW^\chi_F)$.

We will use the following lemma.
\begin{lemma} \label{lem AbramovRohklin}
Let $\chi \in \widehat{\cL}(A, \mu)$, $a \in A$ and let $\cE$ be an $\ta(a)$-invariant $\sigma$-algebra,  then we have
\aryst
 h_\mu(a\mid \cE \vee  \cW^{\chi} ) &\leq& h_\mu(a\mid \cE \vee \cF  \vee  \cW^{\chi}  ) +
h_{m_{G/\Gamma}}(a \mid  \cW^{\chi}_G  ) \\
&\leq&   \sum_{\lambda \in \cL_F(A, \mu) \cap \chi } \lambda(a) \dim E_{F, \lambda} + h_{m_{G/\Gamma}}( a \mid \cW^{\chi}_G).
\earyst
\end{lemma}
\begin{proof}
 The second inequality follows from Ledrappier-Young's entropy bound in \cite{LedYou2}.

The first inequality is  similar to   the coarse Abramov-Rohklin's formula \cite[Theorem 13.6]{BRHW1} (see also \cite[Proposition 3.3]{BRHW}).    We give here a sketch of the proof.

If there exists no $\beta \in \Sigma$ with $\chi = [\beta]$, then we have $\cW^{\chi} = \cW^{\chi} \vee \cF$. The inequality becomes trivial.

Now let us assume that $\chi = [\beta]$ for some  $\beta \in \Sigma$.
Let $\xi$ be an increasing (i.e. with Markov property) measurable partition of $G/\Gamma$ subordinate to the partition into $G_{ \chi }$-orbits.  Let $\eta = \pi^{-1} \xi$ be the associated  increasing measurable partition of $M^\alpha$.  Let $\zeta$ be an  increasing measurable partition of $M^\alpha$ subordinate to $\cW^{\chi }$-manifolds with $\eta\prec \zeta$.
We have
\aryst
h_\mu (a\mid \calE \vee \cW^\chi)
&=& h_\mu(a\mid  \calE\vee\zeta)\\
&=& h_\mu(a\mid  \calE\vee\zeta\vee \eta)\\
&\leq& h_\mu(a\mid  \calE\vee \zeta\vee \eta^\Z)+  h_\mu(a\mid  \calE\vee \eta) \\
&=& h_\mu(a\mid \calE\vee \cW^\chi \vee \cF)+  h_{m_{G/\Gamma}}(a\mid \xi)
\earyst
Here the inequality is standard and can be found, for instance, in \cite[Lem 6.1]{KRH}.  The equivalence of $h_\mu(a\mid  \calE\vee \eta)=h_{m_{G/\Gamma}}(a\mid \xi) $ follows as the Haar measure on $G/\Gamma$ is $a$-ergodic.
\end{proof}

\section{Minimal parabolic subgroups and invariant metrics}

In this section, we will prove the following general theorem, which reduces Conjectures \ref{C:Zimmer1} and \ref{C:Zimmer2} to the study of invariant measures on $M^\alpha$ under a minimal parabolic subgroup of $G$.

\begin{theorem}\label{T:general}
Let $G$ be a connected semisimple Lie group with finite center and without simple factors of real rank less than $2$, $\Gamma\subset G$ be a lattice, $M$ be a compact manifold, and $\alpha:\Gamma\to\Diff^{\reg}(M)$ be an action.  Then the following statements hold:
\begin{itemize}
  \item[(1)] Suppose that $\dim(M)<n(G)$ and
  for some minimal parabolic subgroup $P$ of $G$, every $P$-invariant Borel probability measure on $M^\alpha$ is $G$-invariant. Then $\alpha$ preserves a continuous Riemannian metric.
  \item[(2)] Suppose, in addition to the assumptions in part (1), that $\dim(M)<v(G_{\mathrm{cpt}})$. Then $\alpha(\Gamma)$ is finite.
\end{itemize}
\end{theorem}

\subsection{Proof of Theorem \ref{T:general}}
In this subsection we prove Theorem \ref{T:general}. Let $G$, $P$, $\Gamma$, $M$ and $\alpha$ be as in Theorem \ref{T:general}.  Without loss of generality, we may assume that $P$ is the standard minimal parabolic subgroup. In particular, it contains the maximal
$\RR$-split Cartan subgroup $A$ given in Section \ref{sec preliminaries}.

First, let us recall the following result, which is a combination of some main intermediate theorems in \cite{BFH,BFH-new}.

\begin{theorem}[\cite{BFH,BFH-new}]\label{T:BFH}
Suppose that the action $\alpha:\Gamma\to\Diff^{\reg}(M)$ satisfies the following condition: For every $A$-invariant Borel probability measure $\mu$ on $M^\alpha$
with $\pi_*(\mu)=m_{G/\Gamma}$, one has that $\lambda_+^F(\ba,\mu)=0$ for every $\ba\in\La$. Then $\alpha$ preserves a continuous Riemannian metric on $M$.
\end{theorem}

\begin{proof}
First, the assumption implies that $\alpha$ has uniform subexponential growth of derivatives. In fact, if $\alpha$ fails to have uniform subexponential growth of derivatives, then by \cite[Theorem D]{BFH-new},  there exists a maximal $\R$-split Cartan subgroup  $A$ and an $A$-invariant Borel probability measure $\mu$ on $M^\alpha$ with $\pi_*(\mu)=m_{G/\Gamma}$ and $\lambda_+^F(\ba,\mu)>0$ for some $\ba\in\La$.  Since all maximal $\R$-split Cartan subgroups are conjugate and since $\pi_*(\mu)=m_{G/\Gamma}$, the conclusion holds for any maximal $\R$-split Cartan subgroup $A$, which leads to a contradiction.
See also \cite[Proposition 4.6]{BFH}.

Next, the fact that $\alpha$ has uniform subexponential growth of derivatives implies $\alpha$ preserves a continuous Riemannian metric.
When $\alpha$ takes value in $\Diff^{r}(M)$ for some $r \ge 2$, we may obtain a Riemannian metric of class $C^{r - 1 - \eps}$ for any $\eps > 0$.  This is \cite[Thm.\ 2.4]{BFH} (see also \cite[Thm.\ 1.8]{BFH-new}).
When $\alpha$ takes value in $\Diff^{1+\kappa}(M)$ for $\kappa>0$,  consider first the space $C^{\kappa}(M, S^2T^*M)$ of $\kappa$-H\"older continuous sections into the space of  symmetric $2$-tensors.  A Riemannian metric on $M$ induces a norm on each fiber $ S^2T^*_xM$.
 Given $p\ge 1$ and smooth density $m$ on $M$, we equip $C^{\kappa}(M, S^2T^*M)$ with an appropriate fractional Sobolev norm (defined relative to a choice of partition of unity)
and let $H^{p,\kappa}(M, m, S^2T^*M)$ be the completion of $C^{\kappa}(M, S^2T^*M)$ with respect to this norm.  We make the following observations:
 \begin{enumerate}
 \item $H^{p,\kappa}(M, m, S^2T^*M)$ is in $\mathcal E_{10}$ (for the same reason as in \cite[Lem.\ 7.4]{BFH}.
 \item There is a bounded inclusion $H^{p,\kappa}(M, m, S^2T^*M)$ into $C^{\kappa-\frac{\dim{m}}{p}}(M, S^2T^*M)$  assuming ${\kappa-\frac{\dim{m}}{p}}>0$
 \item If  $\alpha \colon \Gamma\to \Diff^{1+\kappa}(M)$ is an action with uniform subexponential growth of derivatives then the induces action on $H^{p,\kappa}(M, m, S^2T^*M)$ has subexponential growth.
\end{enumerate}
From the above, one may use Strong property (T) exactly as in the proof of  \cite[Thm.\ 2.4]{BFH}  to conclude that $\alpha(\Gamma)$ preserves a Rimannian metric that is at least $C^{\kappa-\epsilon}$.
\end{proof}

To prove statement (1) in Theorem \ref{T:general}, we need to show that the assumptions in statement (1) implies the assumption in Theorem \ref{T:BFH}.
Let us also recall the following implication of Zimmer's cocycle superrigidity theorem.

\begin{lemma}\label{L:GP}
Suppose $\dim(M)<n(G)$. Then for every $G$-invariant Borel probability measure $\mu$ on $M^\alpha$ and every $\ba\in\La$, one has $\lambda_+^F(\ba,\mu)=0$.
\end{lemma}

\begin{proof}
Without loss of generality, we may assume that $\mu$ is $G$-ergodic. Then the assumption $\dim(M)<n(G)$ and Zimmer's cocycle superrigidity theorem imply that the fiberwise derivative cocycle for the $G$-action on $M^\alpha$ is $\mu$-a.e. cohomologous to a cocycle that takes values in the orthogonal group $\OO(\dim(M))$. Hence $\lambda_+^F(\ba,\mu)=0$ for every $\ba\in\La$.
\end{proof}

The following lemma is the main technical result in this section.

\begin{lemma}\label{L:PA}
Suppose for every $P$-invariant Borel probability measure $\mu$ on $M^\alpha$, one has that $\lambda_+^F(\ba,\mu)=0$ for every $\ba\in\La$.
Then, for every $A$-invariant Borel probability measure $\mu$ on $M^\alpha$ with $\pi_*(\mu)=m_{G/\Gamma}$,
one also has that $\lambda_+^F(\ba,\mu)=0$ for every $\ba\in\La$.
\end{lemma}

We postpone the proof of Lemma \ref{L:PA} and first give the proof of Theorem \ref{T:general}.

\begin{proof}[Proof of Theorem \ref{T:general}]
(1) The assumptions and Lemma \ref{L:GP} imply that the assumption of Lemma \ref{L:PA} holds. It then follows from Lemma \ref{L:PA} and Theorem \ref{T:BFH} that $\alpha$ preserves a continuous Riemannian metric.

(2) Same as \cite[Section 2.3]{BFH}.
\end{proof}

To prove Lemma \ref{L:PA}, we first establish some preliminaries. Consider the height function $\htt:\Sigma\to\ZZ$, which is defined as follows: If $\beta=\sum_{\alpha\in\Pi}k_\alpha\alpha$, then
$$\htt(\beta)=\sum_{\alpha\in\Pi}k_\alpha.$$
We enumerate $\Sigma^+$ as $\{\beta_1,\ldots,\beta_m\}$ such that $$\htt(\beta_1)\ge\htt(\beta_2)\ge\cdots\ge\htt(\beta_m),$$
and denote
$$\Lu_0=\{0\}, \qquad \Lu_k=\bigoplus_{i=1}^k\Lg_{\beta_i}, \qquad 1\le k\le m.$$
Note that $\Lu_m=\Ln$.

\begin{lemma}
Each $\Lu_k$ is an ideal of $\Ln$.
\end{lemma}

\begin{proof}
We may assume $k\ge1$. Then
$$[\Lu_k,\Ln]=\sum_{i=1}^k\sum_{j=1}^m[\Lg_{\beta_i},\Lg_{\beta_j}]\subset\sum_{i=1}^k\sum_{j=1}^m\Lg_{\beta_i+\beta_j}.$$
Here we use the convention that $\Lg_\beta=\{0\}$ if $\beta\in\La^*\sm(\Sigma\cup\{0\})$.
Let $i\in\{1,\ldots,k\}$ and $j\in\{1,\ldots,m\}$ be such that $\beta_i+\beta_j\in\Sigma$. Then
$$\htt(\beta_i+\beta_j)=\htt(\beta_i)+\htt(\beta_j)\ge\htt(\beta_k)+1.$$
It follows that $\beta_i+\beta_j\in\{\beta_1,\ldots,\beta_k\}$. Thus $\Lg_{\beta_i+\beta_j}\subset\Lu_k$. This implies that $[\Lu_k,\Ln]\subset\Lu_k$, hence completes the proof.
\end{proof}

For $0\le k\le m$, let $U_k$ denote the connected normal Lie subgroup of $N$ with Lie algebra $\Lu_k$. Note that $U_k$ is also a normal subgroup of $P=CAN$ (see \cite[page 456]{Kn}).

We will also use the following key theorem from \cite{BRHW}.

\begin{theorem}[\cite{BRHW}]\label{T:BRHW}
Let $\mu$ be a Borel probability measure on $M^\alpha$, and let
\ary \label{eq def stab G}
Q = \Stab_G(\mu) := \{ g \in G : g_* \mu = \mu \}.
\eary
\begin{itemize}
  \item[(1)] Suppose $\mu$ is $A$-invariant, $A$-ergodic, and $\pi_*(\mu)=m_{G/\Gamma}$.  Then $\Lg_{[\beta]}\subset\Lie(Q)$ for every $[\beta]\in\widehat\Sigma\sm\widehat\cL_F(A,\mu)$.
  \item[(2)] Suppose $\mu$ is $P$-invariant and $P$-ergodic (thus $Q$ is a standard parabolic subgroup).
  Then $\widehat\Sigma\sm\widehat\Sigma_Q\subset\widehat\cL_F(A,\mu)$, and hence
      $$\dim(M)\ge|\widehat\Sigma\sm\widehat\Sigma_Q|.$$
\end{itemize}
\end{theorem}

\begin{proof}
Item (1) is \cite[Proposition 5.1]{BRHW}. To see that (2) follows from (1), it suffices to notice that the $P$-action on $G/\Gamma$ is uniquely ergodic, thus $\pi_*(\mu)$ is $m_{G/\Gamma}$  (which is $A$-ergodic by Moore's ergodicity theorem), and hence the projection of $\mu$-a.e. $A$-ergodic component of $\mu$ to $G/\Gamma$ is $m_{G/\Gamma}$.
\end{proof}

We are now prepared to prove Lemma \ref{L:PA}.

\begin{proof}[Proof of Lemma \ref{L:PA}]
Assume, to the contrary, that there is an $A$-invariant Borel probability measure $\mu_0$ on $M^\alpha$ such that $\pi_*(\mu_0)=m_{G/\Gamma}$
and $\lambda_+^F(\ba_0,\mu_0)>0$ for some $\ba_0\in \La$. We prove the claim that for $0\le k\le m$, there exists an $AU_k$-invariant Borel probability measure $\mu_k$ on $M^\alpha$ such that $\pi_*(\mu_k)=m_{G/\Gamma}$
and $\lambda_+^F(\ba_k,\mu_k)>0$ for some $\ba_k\in \La$.

For $k=0$, the claim is our assumption. Assume $1\le k\le m$ and the claim holds for $k-1$, namely, there exist an $AU_{k-1}$-invariant Borel probability measure $\mu_{k-1}$ on $M^\alpha$ and $\ba_{k-1}\in \La$ such that $\pi_*(\mu_{k-1})=m_{G/\Gamma}$
and $\lambda_+^F(\ba_{k-1},\mu_{k-1})>0$. Note that $\mu_{k-1}$-a.e. $A$-ergodic component of $\mu_{k-1}$ is $U_{k-1}$-invariant (see \cite[Proof of Claim 5.2]{BRHW}). Therefore, we may choose an $A$-ergodic component $\bar{\mu}_{k-1}$ of $\mu_{k-1}$ satisfying the following conditions:
\begin{itemize}
  \item $\bar{\mu}_{k-1}$ is $AU_{k-1}$-invariant;
  \item $\pi_*(\bar{\mu}_{k-1})=m_{G/\Gamma}$;
  \item $\lambda_+^F(\ba_{k-1},\bar{\mu}_{k-1})>0$.
\end{itemize}
Let $V_k$ be the connected Lie subgroup of $N$ with Lie algebra $\Lg_{\beta_k}\oplus\Lg_{2\beta_k}$, $(F_i)_{i\in\NN}$ be a F{\o}lner
sequence in the group $AV_k$. We have the following.

\begin{lemma} \label{lem weaklimitexist}
 The sequence  $(F_i*\bar{\mu}_{k-1})_{i\in\NN}$ has uniformly exponentially small mass at $\infty$ (in the sense of \cite{BFH-new}), i.e., there is $\kappa_0 > 0$ such that
 \aryst
 \sup_{i}
 \int e^{\tau_0 d(x_0, x)} d(F_i*\bar{\mu}_{k-1})(x) < \infty.
 \earyst
 Consequently, $(F_i*\bar{\mu}_{k-1})_{i\in\NN}$  has a weak-$*$ subsequential limit $\mu_k$ satisfying $\pi_*(  \mu_{k} ) = m_{G/\Gamma}$.
\end{lemma}

\begin{proof}
Since the expression in the integral depends only on the base point, and  for every $i \in \N$, the projected measure $\pi_*( F_i*\bar{\mu}_{k-1} )$ is the Haar measure,  the inequality is clear.
In particular, the collection of measures $F_i*\bar{\mu}_{k-1}$, $i \geq 0$, are uniformly integrable   with respect to the proper function $e^{\tau_0 d(x_0, \cdot)}$.
 Thus the sequence $F_i*\bar{\mu}_{k-1}$ is precompact in the space of Borel probability measures on $M^{\a}$.  This concludes the proof.
\end{proof}

 We prove that $\mu_k$ given by Lemma \ref{lem weaklimitexist} satisfies the requirements of the claim.

It is clear that $\pi_*(\mu_k)=m_{G/\Gamma}$.
Let us show that $\mu_k$ is $AU_k$-invariant. Note that $AU_k$ is generated by $AV_k$ and $U_{k-1}$. Since $\mu_k$ is $AV_k$-invariant, it suffices to show that $\mu_k$ is $U_{k-1}$-invariant. Note that $U_{k-1}$ is normalized by $AV_k$. It follows that for every $h\in AV_k$ and $u\in U_{k-1}$, we have
$$u_*(h_*\bar{\mu}_{k-1})=h_*((h^{-1}uh)_*\bar{\mu}_{k-1})=h_*\bar{\mu}_{k-1}.$$
Thus $h_*\bar{\mu}_{k-1}$ is $U_{k-1}$-invariant for every $h\in AV_k$. This implies that every $F_i*\bar{\mu}_{k-1}$, and hence $\mu_k$, is $U_{k-1}$-invariant. Hence $\mu_k$ is $AU_k$-invariant.

It remains to show that $\lambda_+^F(\ba_k,\mu_k)>0$ for some $\ba_k\in \La$. We divide the proof into two cases.

\textbf{Case 1.} Assume $\cL_F(A,\bar{\mu}_{k-1})\not\subset\RR\beta_k$. Let $\lambda\in\cL_F(A,\bar{\mu}_{k-1})\sm\RR\beta_k$. We choose $\ba_k\in \La$ such that $\lambda(\ba_k)>0$ and $\beta_k(\ba_k)=0$. Then $\exp(\ba_k)$ centralizes $AV_k$.

By  \cite[Lemma 6.5]{BFH-new} and Lemma \ref{lem weaklimitexist},    we deduce that
$$\lambda_+^F(\ba_k,\mu_k)\ge\lambda_+^F(\ba_k,\bar{\mu}_{k-1})\ge\lambda(\ba_k)>0.$$

\textbf{Case 2.} Assume $\cL_F(A,\bar{\mu}_{k-1})\subset\RR\beta_k$. We will prove that $\bar{\mu}_{k-1}$ is invariant under the simple factor $G_k$ of $G$ containing $V_k$. Assuming this, we then have $\mu_k=\bar{\mu}_{k-1}$, and it is enough to take $\ba_k=\ba_{k-1}$. We now prove the $G_k$-invariance of $\bar{\mu}_{k-1}$. Let $Q=\Stab_G(\bar{\mu}_{k-1})$, $\Lq=\Lie(Q)$. By Theorem \ref{T:BRHW}, we have $\Lg_\gamma\subset\Lq$ for every $\gamma\in\Sigma\sm\RR\beta_k$. Consider the spaces
$$\Lb_1=\bigoplus_{\gamma\in\Sigma\cap\RR\beta_k}\Lg_\gamma, \qquad \Lb_2=\bigoplus_{\gamma\in\Sigma\sm\RR\beta_k}\Lg_\gamma.$$
Then the subalgebra $\langle\Lb_2\ra$ of $\Lg$ generated by $\Lb_2$ is contained in $\Lq$.
It suffices to show that the Lie algebra $\Lg_k$ of $G_k$ is contained in $\langle\Lb_2\ra$. Note that $[\Lc\oplus\La\oplus\Lb_1,\Lb_2]\subset\Lb_2$. So $[\Lc\oplus\La\oplus\Lb_1,\langle\Lb_2\ra]\subset\langle\Lb_2\ra$, and hence
$$[\Lg,\langle\Lb_2\ra]=[\Lc\oplus\La\oplus\Lb_1,\langle\Lb_2\ra]+[\Lb_2,\langle\Lb_2\ra]\subset\langle\Lb_2\ra.$$
It follows that $\langle\Lb_2\ra$ is an ideal of $\Lg$. On the other hand, since the real rank of $\Lg_k$ is at least $2$, we have $\Lg_k\cap\Lb_2\ne\{0\}$. So $\Lg_k\subset\langle\Lb_2\ra$. This proves that $\bar{\mu}_{k-1}$ is $G_k$-invariant, hence completes the proof of the claim.

Since $U_m=N$, the measure $\mu_m$ given in the claim is $AN$-invariant. Let $\mu=C*\mu_m$, $\ba=\ba_m$. Since $C$ centralizes $A$ and normalizes $N$, the measure $\mu$ is invariant under $P=CAN$. Moreover, it follows from \cite[Lemma 3.2]{BFH} that $\lambda_+^F(\ba,\mu)\ge\lambda_+^F(\ba_m,\mu_m)>0$. This contradicts the hypothesis of Lemma \ref{L:PA}. Thus the proof is completed.
\end{proof}

\subsection{The quantity $s(G)$}
For the convenience of using Theorem \ref{T:general}, let us make the following definition.

\begin{definition}
For a connected semisimple Lie group $G$ with finite center, let $s(G)$ denote the maximal positive integer such that the following condition holds: For any compact manifold $M$ with $\dim(M)<s(G)$, any lattice $\Gamma\subset G$ and any action $\alpha:\Gamma\to\Diff^{\reg}(M)$, every Borel probability measure on $M^\alpha$ that is invariant under a minimal standard parabolic subgroup $P$ is in fact $G$-invariant.
\end{definition}

As in \cite{BRHW,BFH}, let $r(G)$ denote minimal resonant codimension, namely, the minimal value of $|\widehat\Sigma\sm\widehat\Sigma_Q|$ as $Q$ varies over all proper standard parabolic subgroups of $G$.
In particular, if $G$ is $\RR$-split, then $r(G)=v(G)$.
We have the following bounds for $s(G)$.

\begin{lemma}\label{L:s-ineq}
We always have $r(G)\le s(G)\le v(G)$.
\end{lemma}

\begin{proof}
To prove the first inequality, it suffices to show that if $\dim(M)<r(G)$, then every $P$-invariant Borel probability measure on $M^\alpha$ is $G$-invariant. This follows from Theorem \ref{T:BRHW}, see also \cite[Proof of Theorem 1.6]{BRHW}.

To prove the second inequality, we take a standard parabolic subgroup $Q\subset G$ with $\dim(G/Q)=v(G)$, and let $M=G/Q$. Let a lattice $\Gamma\subset G$ act on $M$ by translations. Then the map
$$\sigma:G/\Gamma\to M^\alpha, \quad g\Gamma\mt(g,g^{-1}Q)\Gamma$$
is a global section and is $Q$-equivariant. It follows that $\sigma_*(m_{G/\Gamma})$ is $Q$-invariant.
On the other hand, the support of $\sigma_*(m_{G/\Gamma})$ is equal to $\sigma(G/\Gamma)$, which is not a $G$-invariant set. So the image measure is not $G$-invariant.
\end{proof}

Let us highlight the following corollary of Theorem \ref{T:general}.

\begin{corollary}\label{C:metric}
Let $G$ be a connected real semisimple Lie group with finite center and without simple factors of real rank less than $2$. Let $\Gamma\subset G$ be a lattice, $M$ be a compact manifold, and $\alpha:\Gamma\to\Diff^{\reg}(M)$ be an action.
\begin{itemize}
  \item[(1)] If $\dim(M)<s(G)$, then $\alpha$ preserves a continuous Riemannian metric.
  \item[(2)] If $\dim(M)<\min\{s(G),v(G_{\mathrm{cpt}})\}$, then $\alpha(\Gamma)$ is finite.
\end{itemize}
\end{corollary}

\begin{proof}
It suffices to notice that $s(G)\le v(G)\le n(G)$.
\end{proof}

In view of Corollary \ref{C:metric}, to prove Conjecture \ref{C:Zimmer1} for a group $G$, it is enough to prove that $s(G)=v(G)$.
In particular, Theorem \ref{T:main-1} follows from the following statements for $s(G)$.

\begin{theorem} \label{thm valuesofs}
\begin{itemize}
  \item[(1)] Suppose that $G$ is either a connected complex semisimple Lie group without simple factors of rank $1$, or the group $E II$. Then $s(G)=v(G)$.
  \item[(2)]  We have
\begin{align}
s(\SL(n,\HH)) & \ge\begin{cases}
4n-4; \\ 8; \\ 4,
\end{cases} && \begin{matrix}
n\ge5; \\  n=4; \\ n=3,
\end{matrix} \label{E:s-SL} \\
s(\SO^+(m,n)) & \ge\begin{cases}
\min\{m+n-2,n(n+3)/2\};\\ 4,
\end{cases} && \begin{matrix}
n\ge3, \ m\ge n+2; \\  n=2, \ m\ge4, \ \ \ \ \ \ \
\end{matrix} \label{E:s-SO} \\
s(\SU(m,n)) & \ge\begin{cases}
\min\{2m+2n-3,n(n+1)\}; \\ 4,
\end{cases} && \begin{matrix}
m\ge n\ge2, \ (m,n)\ne(2,2); \\  (m,n)=(2,2), \qquad \qquad \quad \ \
\end{matrix} \label{E:s-SU} \\
s(\Sp(m,n)) & \ge\min\{4m+4n-7,n(n+1)\}, && m\ge n\ge2, \label{E:s-Sp} \\
s(\SO^*(2n)) & \ge\min\{4n-7,[n^2/4]\}, && n\ge4. \label{E:s-SO*}
\end{align}
\end{itemize}
\end{theorem}

\begin{proof}[Proof of Theorem \ref{T:main-1}]
The cases (i) and (vi) follow directly from Theorem \ref{thm valuesofs} and Corollary \ref{C:metric}. For cases (ii)--(v), it suffices to notice that
\begin{align*}
v(\SL(n,\HH))&=4n-4, && n\ge 3,\\
v(\SO^+(m,n))&=m+n-2, && m>n\ge2,\\
v(\SU(m,n)) & =\begin{cases}
2m+2n-3; \\ 4,
\end{cases} && \begin{matrix}
m\ge n\ge2, \ (m,n)\ne(2,2); \\  (m,n)=(2,2), \qquad \qquad \quad \ \
\end{matrix} \\
v(\SO^*(2n))&=4n-7, && n\ge4,
\end{align*}
and in each case, Theorem \ref{thm valuesofs} gives $s(G)=v(G)$.
\end{proof}

For the group $\Sp(m,n)$, Corollary \ref{C:metric} and Theorem \ref{thm valuesofs} imply the following result.

\begin{theorem}\label{T:Sp}
Let $G=\Sp(m,n)$ with $6\le n\le m\le \frac14(n^2-3n+6)$.
Let $\Gamma\subset G$ be a lattice, $M$ be a compact manifold with
$$\dim(M)<v(G)-2=4m+4n-7,$$
and $\alpha:\Gamma\to\Diff^{\reg}(M)$ be an action. Then $\alpha(\Gamma)$ is finite.
\end{theorem}

\begin{proof}
It suffices to notice that $s(G)\ge v(G)-2=4m+4n-7$ and $v(G_{\mathrm{cpt}})=4m+4n-4$.
\end{proof}

\subsection{Bounding $s(G)$: Two approaches} \label{subsec twoapproaches}

The remaining part of this paper will be devoted to the proofs of Theorem \ref{thm valuesofs} and Theorem \ref{T:main-2}.
In doing this, we use two different approaches.

\smallskip

\noindent{\bf Superrigidity approach}:

The main result we will use is the following theorem, which was first discovered by Zimmer \cite{Zim}.
The version we use can be found in \cite[Theorem 1.4]{FM}.
 \begin{theorem}[Zimmer's cocycle superrigidity theorem] \label{thm ZCS}
Let $G$ be a higher-rank semisimple Lie group, and
let $(S, \mu)$ be a standard probability measure space.
Assume that $G$ acts ergodically on $S$ preserving $\mu$. Let $\alpha :  G \times S \to GL(m, \R)$, with an integer $m \geq 1$, be an integrable Borel cocycle. Then $\alpha$ is cohomogolous to a cocycle $\beta$ where $\beta(g, x) = \pi(g) c(g,x)$. Here $\pi : G \to GL(m, \R)$ is a continuous homomorphism  and $c : G \times S \to C$ is a cocycle taking values in a compact subgroup $C < GL(m, \R)$ centralizing $\pi(G)$.
\end{theorem}

In Section \ref{sec: superrigidity}, we prove a general result, Proposition \ref{P:real-1}, using Theorem \ref{thm ZCS}. Then we will apply it to the group   $\SL(n,\CC)$.

\smallskip

\noindent{\bf  Measure rigidity approach}:

The measure rigidity result we will use is the following.
\begin{theorem} \label{prop:dim2}
Let $P\subset G$ be the minimal standard parabolic subgroup, $\mu$ be a $P$-invariant $P$-ergodic Borel probability measure on $M^\alpha$, and consider the standard parabolic subgroup $Q:=\Stab_G(\mu)$.
Then for every $\beta \in \Sigma \setminus \Sigma_Q$, there is $\chi \in \widehat\cL_{F}(A, \mu)$
such that $\beta \in  \chi $ and $\dim E^{\chi}_F(x) \geq \min \{ \dim \Lg_{[\beta]},    2 \}$ for $\mu$-a.e. $x$.
\end{theorem}

We will prove Theorem \ref{prop:dim2} in Subsections \ref{subsec Extrainvariance} to \ref{subsec finishinglemmaentropy}.
Then in Subsection \ref{subsec anewdimbound}, we will apply it to complex semisimple groups and the groups $E II$, $\SL(n,\HH)$ and $\SO^+(m,n)$.

\smallskip

The measure rigidity approach can be very useful in handling some marginal cases that are not covered by the superrigidity approach. In Section \ref{sec: comination}, we combine the above two approaches to establish the lower bounds for $s(G)$ for the groups $\SU(m,n)$, $\Sp(m,n)$ and $\SO^*(2n)$.

\section{The superrigidity approach} \label{sec: superrigidity}

In this section, we will prove some general estimates related to Theorem \ref{thm valuesofs}
using Zimmer's cocycle superrigidity theorem.
We then illustrate this approach by showing that Theorem \ref{T:main-1} holds for $G=\SL(n,\CC)$ when $n\ge5$.
For simplicity, we assume throughout this section that $G$ is simple.
As in Theorem \ref{T:general}, we let $\Gamma\subset G$ be a lattice, let $M$ be a compact manifold, and let $\alpha:\Gamma\to\Diff^{\reg}(M)$ be an action.

\subsection{The general superrigidity approach}

\begin{proposition}\label{P:real-1}
Let $G$ be a connected simple Lie group with finite center and with real rank at least $2$, and let $\Gamma$, $M$ and $\alpha$  be as above.
Let $P\subset G$ be the minimal standard parabolic subgroup, $\mu$ be a $P$-invariant $P$-ergodic Borel probability measure on $M^\alpha$, and consider the standard parabolic subgroup $Q:=\Stab_G(\mu)$ (see \eqref{eq def stab G}).  Suppose that the set $\Pi_Q\subset\Pi$, viewed as a subgraph of the Dynkin diagram, has a connected component $\Delta$ with $|\Delta|\ge2$. Let $\Lg_\Delta$ be the (simple) subalgebra of $\Lg$ generated by $\{\Lg_\beta:\beta\in\spann(\Delta)\cap\Sigma\}$. Then
$$\dim(M)\ge|\Pi\sm\Delta|n(\Lg_\Delta),$$
where $n(\Lg_\Delta)$ is the minimal dimension of nontrivial real representations of $\Lg_\Delta$.
\end{proposition}

We first prove

\begin{lemma}\label{L:root-span}
Let $\Sigma$ be an irreducible root system (not necessarily reduced) in a Euclidean space $V$, and $W\subset V$ be a nontrivial subspace. Then
$$\spann(\Sigma\sm(W\cup W^\bot))=V.$$
\end{lemma}

\begin{proof}
Suppose to the contrary that $U:=\spann(\Sigma\sm(W\cup W^\bot))\ne V$. Since $\Sigma$ is irreducible, we have $U\ne\{0\}$. By the irreducibility again, we can choose $\alpha\in\Sigma\sm(U\cup U^\bot)$. Note that $\Sigma\subset W\cup W^\bot\cup U$. So $\alpha\in W\cup W^\bot$. Since $\alpha\notin U^\bot$, we can choose $\beta\in\Sigma\sm(W\cup W^\bot)$ such that $\langle\alpha,\beta\rangle<0$. It follows that $\alpha+\beta\in\Sigma$. Let $q\ge1$ be the maximal integer such that $q\alpha+\beta\in\Sigma$. Since $\alpha\notin U$ and $\beta\in U$, we have $q\alpha+\beta\notin U$, and hence $q\alpha+\beta\in W\cup W^\bot$. We claim that $\langle\alpha,q\alpha+\beta\ra=0$. In fact, recall that $\alpha,q\alpha+\beta\in W\cup W^\bot$. If $\alpha\in W$ (resp. $\alpha\in W^\bot$), then it follows from $\beta\notin W\cup W^\bot$ that $q\alpha+\beta\notin W$ (resp. $q\alpha+\beta\notin W^\bot$), and hence $q\alpha+\beta\in W^\bot$ (resp. $q\alpha+\beta\in W$). This verifies the claim. It follows from the claim that $\langle\alpha,\beta\rangle=-q\langle\alpha,\alpha\rangle$. Let $\sigma_\alpha$ denote the orthogonal reflection with respect to the hyperplane $\alpha^\bot$. Then $\sigma_\alpha(\beta)=\beta-\frac{2\langle\alpha,\beta\rangle}{\langle\alpha,\alpha\rangle}\alpha=2q\alpha+\beta\in\Sigma$, which contradicts the maximality of $q$. This proves the lemma.
\end{proof}

\begin{proof}[Proof of Proposition \ref{P:real-1}]
For $\beta\in\Sigma$, let $\bh_\beta\in\La$ be such that $\langle \bh_\beta,\ba\ra=\beta(\ba)$ for all $\ba\in\La$, where $\langle\cdot,\cdot\ra$ is the Killing form of $\Lg$. Then $\RR \bh_\beta\subset[\Lg_\beta,\Lg_{-\beta}]\subset\Lc\oplus\RR \bh_\beta$ (see \cite[(7.63)]{Kn}). Let $\La_\Delta=\spann\{\mathbf{h}_\alpha:\alpha\in\Delta\}$, $\Lc_\Delta=\Lc\cap\Lg_\Delta$. Then $\Lg_\Delta=\Lc_\Delta\oplus\La_\Delta\oplus\bigoplus_{\beta\in\spann(\Delta)\cap\Sigma}\Lg_\beta$. Thus $\La_\Delta$ is a maximal $\RR$-diagonalizable subalgebra of $\Lg$, and the restricted root system $\Sigma(\Lg_\Delta,\La_\Delta)$ is isomorphic to $\spann(\Delta)\cap\Sigma$, for which $\Delta$ is a set of simple roots. Since $\Delta$ is connected, $\Lg_\Delta$ is a simple Lie algebra of real rank $|\Delta|\ge2$. Note that the orthogonal complement $\La_\Delta^\bot$ of $\La_\Delta$ in $\La$ is equal to $\bigcap_{\alpha\in\Delta}\ker(\alpha)$. So $[\La_\Delta^\bot,\Lg_\Delta]=0$.

If $\Delta=\Pi$, the proposition follows immediately from Zimmer's cocycle superrigidity theorem.
So we assume $\Delta\ne\Pi$. Let $W=\spann(\Delta)$. Since $\Delta$ is a connected component of $\Pi_Q$, we have
$$\spann(\Pi_Q)\cap\Sigma\subset\spann(\Delta)\cup\spann(\Pi_Q\sm\Delta)\subset W\cup W^\bot.$$
Thus, it follows from Theorem \ref{T:BRHW} and Lemma \ref{L:root-span} that
\begin{align*}
\spann(\cL_F(A,\mu)\sm W^\bot)&\supset\spann((\Sigma\sm\Sigma_Q)\sm W^\bot)\\
&=\spann((\Sigma\sm(\spann(\Pi_Q)\cap\Sigma))\sm W^\bot)\\
&\supset\spann(\Sigma\sm(W\cup W^\bot))=\La^*.
\end{align*}
In particular, there exist $\lambda_1,\ldots,\lambda_k\in\cL_F(A,\mu)\sm W^\bot$, where $k=|\Pi\sm\Delta|$, such that $\{\lambda_i|_{\La_\Delta^\bot}:1\le i\le k\}$ is a basis of $(\La_\Delta^\bot)^*$. Note that $\lambda_i\notin W^\bot$ means that $\lambda_i|_{\La_\Delta}\ne0$.

Let $G_\Delta$, $A_\Delta$ and $A_\Delta^\bot$ be the connected Lie subgroups of $G$ with Lie algebras $\Lg_\Delta$, $\La_\Delta$ and $\La_\Delta^\bot$, respectively. Consider the   Oseledects splitting $E_F=\bigoplus_{\chi\in  \widehat \cL_F(A_\Delta^\bot,\mu)}E_F^\chi$ for the $A_\Delta^\bot$-action on $E_F$ with respect to $\mu$, where $\cL_F(A_\Delta^\bot,\mu)\subset(\La_\Delta^\bot)^*$ is the corresponding set of Lyapunov functionals, and $\widehat \cL_F(A_\Delta^\bot,\mu)$ is the corresponding coarse Lyapunov functionals.
Then
$\chi_i:= [ \lambda_i|_{\La_\Delta^\bot} ] \in \widehat \cL_F(A_\Delta^\bot,\mu)$ and $E_F^{[\lambda_i]}\subset E_F^{\chi_i}$.
Since $G_\Delta$ commutes with $A_\Delta^\bot$, each $E_F^{\chi_i}$ is $G_\Delta$-invariant. Note that $\lambda_i|_{\La_\Delta}\ne0$ is a Lyapunov functional for the $A_\Delta$-action on $E_F^{\chi_i}$ with respect to $\mu$. By Zimmer's cocycle superrigidity theorem (Theorem \ref{thm ZCS}), we have $\dim E_F^{\chi_i}\ge n(\Lg_\Delta)$. It follows that
$$\dim(M)=\dim E_F\ge\dim\bigoplus_{i=1}^kE_F^{\chi_i}\ge kn(\Lg_\Delta).$$
This proves the proposition.
\end{proof}

\subsection{Application to $\SL(n,\CC)$}
We illustrate the superrigidity approach by showing that Theorem \ref{T:main-1} holds for $\SL(n,\CC)$, $n\ge5$. By Corollary \ref{C:metric} and Lemma \ref{L:s-ineq}, it suffices to prove
\begin{equation}\label{E:superrigidity-SL}
  s(\SL(n,\CC))\ge v(\SL(n,\CC))=2n-2, \qquad n\ge5.
\end{equation}

\begin{proof}[Proof of \eqref{E:superrigidity-SL}]
Let $G=\SL(n,\CC)$, $n\ge5$. We can take $\La$ to be the space of real diagonal matrices of trace zero. Then
$$\Sigma=\{e_i-e_j:1\le i\ne j\le n\},$$
where $e_i\in\La^*$ is the functional sending a diagonal matrix to its $i$-th diagonal entry.
The set of simple roots $\Pi=\{\alpha_1,\ldots,\alpha_{n-1}\}$ can be chosen as
$$\alpha_1=e_1-e_2, \quad \ldots, \quad \alpha_{n-1}=e_{n-1}-e_n.$$
We need to prove that if $\alpha$ is a $C^{\reg}$ action of a lattice $\Gamma\subset G$ on a compact manifold $M$ such that there is a $P$-invariant Borel probability measure on $M^\alpha$ that is not $G$-invariant, then
$$\dim(M)\ge2n-2.$$
Without loss of generality, we may assume that $\mu$ is $P$-ergodic. Let $Q=\Stab_G(\mu)$. Then $\Pi_Q\ne\Pi$. There are 2 cases.

\textbf{Case 1.} Suppose that $\Pi_Q$ has a connected component $\Delta$ with $|\Delta|\ge2$, say $\Delta=\{\alpha_p,\alpha_{p+1},\ldots,\alpha_q\}$, where $1\le p<q\le n-1$ and $(p,q)\ne(1,n-1)$. Then $\Lg_\Delta\cong\Lsl(q-p+2,\CC)$. By Proposition \ref{P:real-1}, we have
$$\dim(M)\ge |\Pi\sm\Delta|n(\Lg_\Delta)=2(n-q+p-2)(q-p+2)\ge 2n-2.$$

\textbf{Case 2.} Suppose every connected component of $\Pi_Q$ consists of a single root. Then $\Sigma\sm\Sigma_Q=-(\Sigma^+\sm\Pi_Q)$ and $|\Pi_Q|\le[n/2]$. By Theorem \ref{T:BRHW}, we have
$$\dim(M)\ge |\Sigma\sm\Sigma_Q|=|\Sigma^+\sm\Pi_Q|\ge n(n-1)/2-[n/2]\ge 2n-2.$$
This proves \eqref{E:superrigidity-SL}.
\end{proof}

\section{The measure rigidity approach} \label{sec: measurerigidity}

From this section, we will approach Theorems \ref{thm valuesofs} and \ref{T:main-2} using measure rigidity method.
We will deduce Item (1), \eqref{E:s-SL} and \eqref{E:s-SO} in   Theorem \ref{thm valuesofs} and prove Theorem \ref{T:main-2}.

We let $G, \Gamma, M, \alpha$ be as in Theorem \ref{T:general}.

We will prove Theorem \ref{prop:dim2} in Subsection \ref{subsec Extrainvariance} to \ref{subsec finishinglemmaentropy}.  Then we will give the applications of Theorem \ref{prop:dim2} in Subsection \ref{subsec anewdimbound}.

\subsection{An entropy equation}

 Let $\beta$ be as in Theorem \ref{prop:dim2}.
 In this subsection, we make no assumptions on the dimension of $E^{[\beta]}_F$.

Let $A'_{\beta} = \exp(\ker(\beta))$.
Let $\cE_{\beta}$ denote the $A'_{\beta}$-ergodic decomposition of $\mu$.
By definition,  $\mu^{\cE_{\beta}}_x$ is the conditional measure of the $A'_{\beta}$-ergodic component of $\mu$ at $x$.

 Let $\cF$ denote the partition of $M^\alpha$ into fibers of $\pi$.
Observe that both the ergodic decomposition ${\cE_\beta}$ and $\cF$ are  $A$-invariant.  Given $a\in A$, we may define the conditional entropy $$h_\mu(a\mid {\cE_\beta})$$ relative to the ergodic decomposition  as well as the fiberwise conditional entropy $$h_\mu(a\mid {\cE_\beta}\vee \cF)$$ relative to the ergodic decomposition.
Also, given  $a \in A$ with $\beta(a)>0$ (recall the notation in Subsection \ref{subsec Structure of semisimple Lie groups}), we may define the entropy contributions of $[\beta]$ to be
$$h_\mu( a \mid {\cE_\beta} \vee  \cW^{[\beta]})$$
and
$$h_\mu( a \mid {\cE_\beta}\vee \cF\vee \cW^{[\beta]})$$
along coarse Lyapunov manifolds and fiberwise coarse Lyapunov manifolds.

Fix and $a  \in A$ with $\beta(a)>0$.
We have the bound
\begin{equation}\label{basicLYbound}0\le h_\mu(a\mid {\cE_\beta} \vee  \cW^{[\beta]})- h_\mu(a\mid {\cE_\beta}\vee \cF\vee \cW^{[\beta]})\le \beta( a) ( \dim \Lg_\beta + 2	 \dim \Lg_{2\beta} ). \end{equation}
Indeed, the lower bound \eqref{basicLYbound} follows from the coarse Lyapunov Abramov-Rohlin  formula \cite[Theorem 13.6]{BRHW1}  and the upper bound follows from  Lemma \ref{lem AbramovRohklin}.
We assert the upper bound in \eqref{basicLYbound} is, in fact, an equality.

\begin{lemma} \label{lemma:entropy}
	Let $\mu$, $P$, $Q$ be as in Theorem \ref{prop:dim2}. Let $\gamma \in \Sigma \setminus \Sigma_Q$ such that $\frac 1 2  \gamma \notin \Sigma$.  Fix $a \in A$ with $\gamma(a)>0$.
	Then
\begin{equation}\label{eq:entgaplower}
	\gamma( a) ( \dim \Lg_\gamma + 2	 \dim \Lg_{2 \gamma } ) =
	h_\mu(a\mid {\cE_\gamma} \vee  \cW^{[ \gamma ]})-
	h_\mu(a\mid {\cE_\gamma }\vee \cF \vee \cW^{[ \gamma ]}).
\end{equation}
\end{lemma}
We defer the proof of Lemma \ref{lemma:entropy} to Subsection \ref{subsec finishinglemmaentropy}.

The following corollary will be used in Section \ref{sec: comination}.
\begin{corollary} \label{cor entropyinequalitytononatomic}
Let $\gamma$ be as in Lemma \ref{lemma:entropy}.
Suppose there is $a \in A$ with $\gamma(a) > 0$ and $\kappa>0$ such that $\lambda(a)=\kappa$ for all $\lambda \in   \cL_{F}\cap[\gamma]$.  If
\aryst
\sum_{\lambda \in   \cL_{F}(A, \mu) \cap [\gamma]} \lambda(a) \dim E_{F, \lambda} <  \gamma( a) ( \dim \Lg_\gamma + 2	 \dim \Lg_{2 \gamma } ),
\earyst
then $\mu^{G_{[\gamma]}}_x$ is non-atomic for $\mu$-a.e. $x$.
\end{corollary}

\begin{proof}
Assume to the contrary that $\mu^{G_{[\gamma]}}_x$ is atomic for $\mu$-a.e. $x$.  Then the entropy of $a$ along $G_{[\gamma]}$-orbits vanishes, $h_{\mu}(a \mid G_{[\gamma]})=0$.
Moreover, restricted to  $A'_{\beta}$-ergodic components, the entropy similarly vanishes, $h_{\mu}(a \mid {\cE_\gamma} \vee G_{[\gamma]})=0$.

The subfoliation of  $\cW^{[\gamma]}$-leaves in $M^\alpha$ by $G_{[\gamma]}$-orbits is a smooth foliation.  Since there is a single Lyapunov exponent $\kappa>0$ in $\cW^{[\gamma]}$-leaves transverse to the Lipschitz foliation by $G_{[\gamma]}$-orbits, one may adapt the proof of \cite[Prop.\ 5.1]{LedYou1} to bound the difference  $h_{\mu}( a \mid {\cE_\gamma} \vee\cW^{[\gamma]})-h_{\mu}(a \mid {\cE_\gamma} \vee G_{[\gamma]})$ by $$h_{\mu}( a \mid {\cE_\gamma} \vee\cW^{[\gamma]})-h_{\mu}(a \mid {\cE_\gamma} \vee G_{[\gamma]}) \le \kappa \sum_{\lambda \in   \cL_{F}(A, \mu) \cap [\gamma]}  \dim E_{F, \lambda}.$$
However this contradicts Lemma \ref{lemma:entropy}.
Indeed, we assumed $h_{\mu}(a \mid {\cE_\gamma} \vee G_{[\gamma]}) =0$ and so
\begin{align*}
h_{\mu}( a \mid {\cE_\gamma} \vee \cW^{[\gamma]})
&\le \kappa \sum_{\lambda \in   \cL_{F}(A, \mu) \cap [\gamma]}  \dim E_{F, \lambda}\\
&=   \sum_{\lambda \in   \cL_{F}(A, \mu) \cap [\gamma]} \lambda(a)  \dim E_{F, \lambda}\\
&<  \gamma( a) ( \dim \Lg_\gamma + 2	 \dim \Lg_{2 \gamma } )\\
&\le h_\mu(a\mid {\cE_\gamma} \vee  \cW^{[ \gamma ]}).
\end{align*}
This is a contradiction.
Thus $\mu^{G_{[\gamma]}}_x$ is non-atomic for $\mu$-a.e. $x$.
\end{proof}

\subsection{Extra invariance from homogenity} \label{subsec Extrainvariance}

Let $\beta$ be as in Theorem \ref{prop:dim2}. In particuliar, $\mu$ is not $G_{[\beta]}$-invariant.

Let $\scrV$ denote an $A$-invariant measurable foliation that is (essentially) saturated by $G_{[\beta]}$-orbits.
For any $x$ such that $\mu^{\cE_\beta}_x$ is defined, we may define the leafwise measure $(\mu^{\cE_\beta}_x)^{\scrV}_y$ for $\mu^{\cE_\beta}_x$-a.e. $y$. Since we have $\mu^{\cE_\beta}_x = \mu^{\cE_\beta}_y$ for $\mu$-a.e. $x$ and $\mu^{\cE_\beta}_x$-a.e. $y$, we see that  for $\mu$-a.e. $y$ the measure $(\mu^{\cE_\beta}_y)^{\scrV}_y = (\mu^{\cE_\beta}_x)^{\scrV}_y$ is defined, which will be abbreviated as $\mu^{\scrV, {\cE_\beta}}_y$. We see that $\mu^{\scrV, {\cE_\beta}}_x$ is the leafwise measure of $\mu^{\cE_{\beta}}_x$ along $\scrV$ at $x$ (recall Section \ref{sec conditionalmeasure}) whenever $\mu^{\scrV, {\cE_\beta}}_x$ is defined.
In practice, we will only consider the following:
\enmt
\item $\scrV(x)=  \cW_{G_{[\beta]}}(x) = G_{[\beta]}\cdot x$ is the $G_{[\beta]}$-orbit of $x$ for $\mu$-a.e. $x$. In this case, we will also abbreviate $\mu^{\scrV, {\cE_\beta}}_x$ as $\mu^{G_{[\beta]},{\cE_\beta}}_x$,
\item $\scrV(x)= \cW^{[\beta]}( x)$ is the full $[\beta]$-coarse Lyapunov manifold through $x$ for $\mu$-a.e. $x$.
\eenmt

\begin{proposition}\label{prop:moreentr}
Suppose for $\mu$-a.e.\ $x$,  the leafwise measure
$\mu^{\scrV,{\cE_\beta}}_x$ is in the Lebesgue class on an embedded submanifold $N_x\subset \scrV(x)$  with respect to the leafwise topology on $\scrV(x)$.    Moreover, we assume that for $\mu$-a.e.\ $x$ that  $T_x(G_{[\beta]}\cdot x)\cap T_xN_x\neq\{0\}$.   Then the measure $\mu$ is  $G_{[\beta]}$-invariant.
\end{proposition}

\subsection{Proof of Proposition \ref{prop:moreentr}}
Proposition \ref{prop:moreentr}  follows directly from the following two lemmas, Lemma \ref{lem:entrop} and  Lemma \ref{lem:alge}, below.

In the first lemma, we canonically identify the tangent space $T_x( G_{[\beta]} \cdot x)$ of the $G_{[\beta]}$-orbit of $x$ in $\cW^{[\beta]}(x)$ with the Lie algebra $\Lg_{[\beta]}$ of $G_{[\beta]}$.  

\begin{lemma}\label{lem:entrop}
Suppose for $\mu$-a.e.\ $x$, the measure $\mu^{G_{[\beta]},{\cE_\beta}}_x$ is in the Lebesgue class on an embedded submanifold $N_x\subset \cW^{[\beta]}(x)$ with respect to the leafwise topology on $ \cW^{[\beta]}(x)$.
Then for $\mu$-a.e.\ $x$ the following holds:

Given a non-zero $Y\in T_x(G_{[\beta]}\cdot x)\cap T_xN_x$, write $Y = Y_\beta+ Y_{2\beta}$ where $Y_\beta\in \Lg_\beta$ and  $Y_{2\beta}\in \Lg_{2\beta}$.
Let $\gamma = \beta$ or $\gamma= 2\beta$ and let $Y' =Y_\gamma$.  Then for $\mu$-a.e.\ $y$, the  $A'_\beta$-ergodic component  $\mu^{{\cE_\beta}}_y$ is invariant under the $1$-parameter subgroup $H=\exp(tY')\subset G_{[\beta]}$.
\end{lemma}

The second fact follows from the structure of $\mathfrak{sl}(2)$-triples; see  \cite[Lemma 7.73]{Kn}.

\begin{lemma}\label{lem:alge}
Let $\beta\in \Sigma$ and let $\Lh$ be a Lie subalgebra such that  $\Lg_{[-\beta]}\subset \Lh $  and  $Y\in \Lh$ for some non-zero $Y\in \Lg_{[\beta]}$.  Then  $\Lg_{[\beta]}\subset \Lh$.
\end{lemma}

\begin{proof}[Proof of Lemma \ref{lem:entrop}]
We first note that for $\mu$-almost every $x$, by Oseledects theorem we have $Y_\beta, Y_{2\beta}\in T_x(G_{[\beta]}\cdot x)\cap T_xN_x$.  Fix a choice of $\gamma$ and set $Y= Y_\gamma$.

Observe that the $\R$-span of $Y$ is contained in $T_x(G_{[\beta]}\cdot x)\cap T_xN_x$.
Given any $a\in A$, we have $\Ad(a) Y= \gamma(a) Y$.  In particular, as $A$ permutes $A'_\beta$-ergodic components, we have that   $Y\in T_y(G_{[\beta]}\cdot y)\cap T_yN_y$ for $\mu$-a.e.\ $y$.  Since $N_y$ is a submanifold and $\mu_y^{\cE'_\beta}$ has   dense support  in $N_y$, we may pass to the closure and conclude  for $\mu$-almost every $x$, every choice of such $\gamma$ and $Y$ at $x$, $\mu$-almost every $y$, and every $z\in N_y$, that $Y\in T_z(G_{[\beta]}\cdot z) \cap T_zN_y$ and thus the orbit  $H\cdot z$ is contained in $N_y$ for every $z\in N_y$.

Let $\cH$ denote the partition into $H$-orbits.
By assumption and Fubini's theorem, the leafwise measure $\mu^{\cH,\cE_\beta}_y$ is absolutely continuous for $\mu^{\cE_\beta}_x$-a.e.\ $y$.
Consider any  $a_1\in A\sm A'_\beta$ with $\beta(a_1)>0$. 	 We have $h_\mu(a_1\mid \cH\vee {\cE_\beta}) = \gamma(a_1)$.  It follows from the explicit density function in \cite{LedYou1} that  $\mu^{\cH,\cE_\beta}_y$ 
is the Haar measure along the parameterized $H$-orbit for $\mu_x^{\cE_{\beta}}$-a.e. $y$.  The result then follows.
\end{proof}

\begin{proof}[Proof of Proposition \ref{prop:moreentr}]
If $\scrV = \cW^{[\beta]}$ as in (2), then the hypothesis of the proposition implies that the leafwise measure $\mu^{G_{[\beta]}, {\cE_\beta}}_x$ is in the Lebesgue class on an embedded submanifold $N_x\subset  \cW^{[\beta]}(x)$. Since $T_x(G_{[\beta]}\cdot x)\cap T_xN_x\neq\{0\}$, we may take a non-zero $Y \in T_x(G_{[\beta]}\cdot x)\cap T_xN_x$ and apply Lemma \ref{lem:entrop} to such $Y$.
We deduce that for $\mu$-a.e. $y$, $\mu^{{\cE_\beta}}_y$ is invariant under the $1$-parameter subgroup $H=\exp(tY')\subset G_{[\beta]}$ for some $Y' \in \Lg_{[\beta]}$.
Then we have $\Lg_{[\beta]} \cap \Lq \neq \{0\}$ where $\Lq$ is the Lie algebra of $Q$. Since $P \subset Q$ by hypothesis, we have $\Lg_{[-\beta]} \subset \Lq$. By Lemma \ref{lem:alge}, we deduce that $\Lg_{[\beta]} \subset \Lq$ as well. This implies that $\mu$ is $G_{[\beta]}$-invariant.
\end{proof}

\subsection{Measure rigidity and homogenity} \label{subsec measurerigidityhomogenity}

Let $\mu$ be a $P$-invariant, $P$-ergodic Borel probability measure on $M^\alpha$ as in Theorem \ref{prop:dim2}.
As before, we let $\beta$ be a restricted root in $\Sigma \setminus \Sigma_Q$, i.e., $\mu$ is not $G_{[\beta]}$-invariant. Without loss of generality, let us assume that $\beta / 2 \notin \Sigma$.
Moreover, throughout this subsection, we assume that $\dim E^{[\beta]}_F(x) = 1$ for $\mu$-a.e. $x$.

Let $\scrV$ be as in Subsection \ref{subsec Extrainvariance} ($\scrV(x) =  G_{[\beta]}\cdot x$
or $\scrV(x)=  \cW^{[\beta]}(x)$).

We have the following.
\begin{lemma} \label{lem measurerigidity}
Assume that we are in one of the following cases:
\enmt
\item $\scrV =  \cW^{[\beta]}_G$, and for $\mu$-a.e. $x$, $\mu^{\scrV, \cE_{\beta}}_x = \mu^{G_{[\beta]}, \cE_{\beta}}_x$ is non-atomic,
\item $\scrV =  \cW^{[\beta]}$, and for $\mu$-a.e. $x$,   $\mu^{G_{[\beta]}, \cE_{\beta}}_x$ is atomic and $\mu^{\scrV, \cE_{\beta}}_x$ is non-atomic.
 \eenmt
Then $\mu^{\scrV,{\cE_\beta}}_x$ is in the Lebesgue class of an embedded submanifold $N_x \subset \scrV(x)$ for $\mu$-a.e. $x$.
\end{lemma}
\begin{proof}
We will only detail the case where $[\beta] \cap \Sigma = \{ \beta, 2 \beta \}$,   as the other case ($[\beta] \cap \Sigma = \{ \beta \}$) follows from a similar but easier argument.

We first detail the proof for (1).

Assume for a positive measure subset of $x\in M^\alpha$ that the leafwise measure  $\mu^{G_{[\beta]},\calE_\beta }_x$ is not a single atom.  By equivariance of $\mu^{G_{[\beta]},\calE_\beta }_x$ under the $A$-action, this holds for $\mu$-a.e.\ $x$.

Fix $\mu'$, an ergodic component of $\mu$ projecting to the Haar measure on $G/\Gamma$.
We have $ \mu^{G_{[\beta]},\calE_\beta }_x=(\mu')^{G_{[\beta]} }_x$.
By \cite[Thm.\ 1]{PS}, we may select $a_0\in A_\beta'$ acting ergodically on $(M^\alpha, \mu')$.

 Note that $a_0$ acts by isometries on $G_{[\beta]}$-orbits and preserves the canonical parametrization by left-translation on orbits: \begin{equation}\label{eq:param}a_0 \cdot (u\cdot x) = u\cdot (a_0\cdot x)\end{equation}
for all $x\in M^\alpha$ and $u\in G_{[\beta]}$.
By the standard measure rigidity arguments, for $\mu'$-a.e.\ $x$ and  $\mu^{G_{[\beta]},\calE_\beta }_x$-a.e.\ $y \in  G_{[\beta]}\cdot x$, there is an isometry $\phi\colon G_{[\beta]}\cdot x\to G_{[\beta]}\cdot x$ such that $\phi(x) = y$ and $\phi_* \left(\mu^{G_{[\beta]},\calE_\beta }_x\right)$ coincides with $\mu^{G_{[\beta]},\calE_\beta }_x$ up to normalization.
In fact, $\phi$ is of the form of   right-translation: if $y= u_y\cdot x$ then for all $u\in G_{[\beta]}$,
$$\phi(u\cdot x) = u\cdot y= u\cdot (u_y\cdot x).$$

It now follows from standard arguments (see for example Lemmas 6.1--6.3 in \cite{KS} following Section 5 of \cite{KS2}) that for $\mu'$-a.e.\ $x$, the measure $\mu^{G_{[\beta]},\calE_\beta }_x$ is the image of the Haar measure on the orbit $V_x\cdot x$ of a connected subgroup $V_x\subset  G_{[\beta]}$.
(Briefly, after slight reformulation below shows the measure $\mu^{G_{[\beta]},\calE_\beta }_x$ is in the Lebesgue class on a connected embedded submanifold $N_x\subset G_{[\beta]}\cdot x$; \eqref{eq:param} shows $T_yN_x\subset \Lg_{[\beta]}$ is constant as $y$ varies in $N_x$.  It follows that $N_x$ is the orbit of a subgroup $V_x$ in $G_{[\beta]}$ and standard entropy arguments show $\mu^{G_{[\beta]},\calE_\beta }_x$ measure is the Haar measure on this $V_x$-orbit.)

Finally, $A$-equivariance of $\mu^{G_{[\beta]},\calE_\beta }_x$ implies the map $x\mapsto \Lie(V_x)$ is $A$-invariant, hence $\mu$-a.e.  constant. It then follows that $\mu$ is invariant under a subgroup $V\subset G_{[\beta]}$.
This concludes the proof of (1).

 Now we proceed to the proof of (2).

Recall that $\scrV =  \cW^{[\beta]}$, and for $\mu$-a.e. $x$, $\mu^{\scrV, \cE_{\beta}}_x$ is non-atomic and $\mu^{G_{[\beta]}, \cE_{\beta}}_x$ is atomic.

We denote by $\cG_1$ the 1-dimensional Lie group of diffeomorphisms of $G_{[\beta]} \times \R$ of form
\aryst
h(u, t) = (u, \lambda_0 t), \quad \lambda_0 \in \R^*.
\earyst
Denote by $\cG$ the finite dimensional Lie group of diffeomorphisms of $G_{[\beta]} \times \R$ of form
\aryst
h(u, t) = (u \cdot u_0, \lambda_0 t + t_0)
\earyst
for some $u_0 \in G_{[\beta]}$, $\lambda_0 \in \R^*$ and $t_0 \in \R$.

Since $\dim E^{[\beta]}_F = 1$, we may construct, using the non-stationary normal form along $\cW^{[\beta]}_F$, a measurable family of $C^{\reg}$ diffeomorphisms $\{  H_x : \cW^{[\beta]}(x) \to G_{[\beta]} \times \R  \}_{x \in M^{\alpha}}$ such that
\enmt
\item
$H_x(x) = (e, 0)$  where $e$ denotes the identity of $G$. Moreover,
$DH_x|_{E^{[\beta]}_F(x)}$ is an isometry,
and  $H_x( \ta(g)(x) ) =  (g,  0)$ for every $g \in G$,
\item for every $a \in \textcolor{teal}{A'_\beta}$ and $\mu$-a.e. $x$,  $H_{\tilde\alpha(a)(x)} \tilde\alpha(a) H_x^{-1} \in \cG_1$,
\item for $\mu$-a.e. $x$, for $\mu^{\cW^{[\beta]}, \cE_{\beta}}_{x}$-a.e. $y$, $H_{y} H_x^{-1} \in \cG$.
\eenmt

For $\mu$-a.e. $x$, we define a Radon measure on $G_{[\beta]} \times \R$ by formula
\aryst
\Psi_x = [ (H_x)_*(\mu^{\cW^{[\beta]}, \cE_{\beta}}_x) ].
\earyst

We have the following lemma.
\begin{lemma} \label{lem returnderivativebound}
For every $\delta > 0$, there exist $C_{\delta} > 1$ and a subset $K \subset M^{\a}$ with $\mu(K) > 1- \delta$ such that for every $x \in K$, every $a \in A'_{\beta}$ with $\tilde\alpha(a)(x) \in K$, we have
\aryst
C_{\delta}^{-1} \leq  \| D \tilde\alpha(a)|_{E^{[\beta]}_F}  \| \leq C_{\delta}.
\earyst
\end{lemma}
\begin{proof}

We fix some $\delta > 0$.

In the following,  we equip $G_{[\beta]}$ with a right-invariant metric.
We let $D_r$ denote the ball centered at $e$ in $G_{[\beta]}$ with radius $r > 0$ with respect to such metric.

By Lemma \ref{lemma:entropy}, we deduce that  $\mu^{\cW^{[\beta]}, \cE_{\beta}}_x$ is not supported on $\cW^{[\beta]}_F(x)$ for $\mu$-a.e. $x$, for otherwise, we would have $0$ on the right hand side of \eqref{eq:entgaplower} for any choices of $a$.
Consequently, for $\mu$-a.e. $x$, $\Psi_x$ is not supported on $\{e\} \times \R$.
Then, using the fact that the $A$-action permutes the $A'_{\beta}$-ergodic components of $\mu$, we deduce
that there exist $R > 0$ and a subset $K_0 \subset M^{\a}$ with $\mu(K_0) > 1 - \delta/2$ such that  for every $x \in K_0$:
\enmt
\item $\mu^{\cW^{[\beta]}, \cE_{\beta}}_x$ and  $H_x$ are defined,
\item the restriction of $\Psi_x$ to $D_1 \times [-R, R]$ is not supported on $\{e\} \times [- R, R]$,
\item $[ \tilde\alpha(a)_* \mu^{\cW^{[\beta]}, \cE_{\beta}}_x] = [ \mu^{\cW^{[\beta]}, \cE_{\beta}}_{\tilde\alpha(a)(x)}]$ for every $x \in K_0$ and every $a \in A'_{\beta}$ with $\ta(a)(x) \in K_0$.
\eenmt

By letting $\theta$ be sufficiently small depending on $\delta$,
we may find a constant $\sigma > 0$ and a subset $K_1 \subset K_0$ with $\mu(K_1) > 1 - 3\delta / 4$ such that
for every $x \in K_1$, we have
\ary \label{eq psixthetapsix1upper}
\frac{\Psi_x(D_{\theta} \times [- R, R])}{\Psi_x(D \times [- R, R])} \leq 1 - \sigma.
\eary

By hypothesis, for $\mu$-a.e. $x$, $\mu^{G_{[\beta]}, \cE_{\beta}}_x$ is a single atom $\delta_x$.
Then there is a compact subset $K_2 \subset K_1$  with $\mu(K_2) > 1- 7\delta/8$ such that for every $x \in K_2$:
\enmt
\item $\mu^{G_{[\beta]}, \cE_{ \beta}}_{x}$ is defined,
\item $\mu^{G_{[\beta]},  \cE_{\beta}}_{x} = \delta_x$.
\eenmt
Given $x \in K_2$, we let $K_x = H_x^{-1}(K_2)$. For every $t \in \R$, the set $K_x \cap (G_{[\beta]} \times \{t\})$ is either empty or a singleton. By compactness of $K_2$, we deduce that for every $x \in K_2$, for all sufficiently small $\kappa > 0$, we have
\ary  \label{eq d1todtheta}
K_y \cap (D_1 \times [-\kappa, \kappa]) = K_y \cap (D_{\theta} \times [-\kappa, \kappa]).
\eary

Let us fix some $x \in K_2$.
The natural projection $p : D_2 \times \R \to \R$ is proper and is measurably one-to-one. The restriction of $\Psi_x$ to $D_2 \times \R$ projects under $p$ to a Radon measure on $\R$, denoted by $\psi_x$. Then almost every $t \in \R$ with respect to $\psi_x$ is a density point with respect to $\psi_x$.
 Moreover, since $p$ is a measurable isomorphism, the map $p^{-1}$ is defined $\psi_x$-almost everywhere.  Consequently, restricted to any interval, $p^{-1}$ is continuous on a subset up to arbitrarily small $\psi_x$-measure.
We deduce that for $\Psi_x$-almost every $(u, t) \in K_x \cap ( D_1 \times \R )$, we have
\ary \label{eq ratiolimitto1}
\lim_{\kappa \to 0} \frac{\Psi_x(K_x \cap ( D_1 \cdot u \times [t-\kappa, t+\kappa] ) )}{\Psi_x( D_1 \cdot u \times [t-\kappa, t+ \kappa] )} = 1.
\eary

  For every $\kappa > 0$ and every $y \in K_2$, we set
\aryst
 c_{\kappa}(y) = \frac{\Psi_y( K_{y, \kappa} ) }{\Psi_y( D_1 \times [-\kappa, \kappa] )} \ \mbox{ where } \
K_{y, \kappa} = K_y \cap (D_1 \times [-\kappa, \kappa]).
\earyst
Note that  $H_x H_y^{-1}$ maps $K_y$ and $D_1 \times [-\kappa, \kappa]$ respectively to $K_x$ and
\aryst
D_1 \cdot u \times [t- a\kappa, t+ a\kappa]
\earyst
for some $u \in G_{[\beta]}$, $t \in \R$ and $a > 0$.  By \eqref{eq ratiolimitto1}, we see that for $\mu^{\cW^{[\beta]}, \cE_{\beta}}_x$-a.e. $y$ in $K_2$,
\ary \label{eq ckappayto1}
\lim_{\kappa \to 0} c_{\kappa}(y) = 1.
\eary

We now select the set $K$ in the proposition.
By \eqref{eq d1todtheta}, \eqref{eq ckappayto1}, and by letting $\kappa > 0$ be sufficiently small, there is a subset $K \subset K_2$ with $\mu(K) > 1 - \delta$  such that for every $y \in K$, we have:
\enmt
\item for every $0 < \kappa' < \kappa$,
\ary \label{eq psixthetapsix1lower}
c_{\kappa'}(y) \geq 1 - \sigma/2,
\eary
\item $K_{y, \kappa'} \subset D_{\theta} \times [-\kappa', \kappa']$  for every $0 < \kappa' < \kappa$.
\eenmt

Now suppose $x \in K$ and $y = \ta(a)(x) \in K$ for some $a \in A'_{\beta}$. Let
\aryst
L = \|  D \tilde\alpha(a)|_{E^{[\beta]}_F} \|.
\earyst
We have
\aryst
H_{\tilde\alpha(a)(x)} \tilde\alpha(a) H_x^{-1}(u, t) = (u, \pm L t).
\earyst
Since $[ \tilde\alpha(a)_* \mu^{\cW^{[\beta]}, \cE_{\beta}}_x] = [ \mu^{\cW^{[\beta]}, \cE_{\beta}}_{y}]$, we have
\aryst
\frac{\Psi_x( D_{\theta} \times [- R, R] ) }{\Psi_x( D_{1} \times [- R, R]  )} = \frac{\Psi_y(D_{\theta} \times [- LR, LR]  ) }{\Psi_y( D_{1} \times [- LR, LR]   )}
\earyst
and
\aryst
\frac{\Psi_x( D_{\theta} \times [- L^{-1}R,   L^{-1}R] ) }{\Psi_x( D_{1} \times [-  L^{-1}R,  L^{-1}R]  )} = \frac{\Psi_y(D_{\theta} \times [- R, R]  ) }{\Psi_y( D_{1} \times [- R, R]   )}.
\earyst
Combined with \eqref{eq psixthetapsix1upper} and \eqref{eq psixthetapsix1lower}, we deduce that
\aryst
\kappa R^{-1} \leq L \leq \kappa^{-1}R.
\earyst
This concludes the proof.
\end{proof}

\smallskip

 We let $K$ be a compact set satisfying  Lemma \ref{lem returnderivativebound} and such that coordinates $x\mapsto H_x $ and measures $x\mapsto \Psi_x$ vary continuously for $x\in K$.   Fix an $A'_{\beta}$-ergodic component $\mu^{\cE_{\beta}}_z$ of $\mu$ and $a \in A'_{\beta}$ such that $\mu^{\cE_{\beta}}_z$ is $a$-ergodic.
For $\mu^{\cE_{\beta}}_z$-a.e.\ $x\in K$ and $\mu^{\cW^{[\beta]}, {\cE_\beta}}_x$-a.e.\ $x'\in
\cW^{[\beta]}(x)\cap K$ there are $y\in K$ and $n_j\to \infty$ and $m_j\to \infty$ with
\begin{enumerate}
\item $ \ta(a^{n_j})(y) \in K$,
\item $ \ta(a^{m_j})(y) \in K$,
\item $ \ta(a^{n_j})(y) \to x$, and
\item $ \ta(a^{m_j})(y) \to x'$.
\end{enumerate} 
 By Lemma \ref{lem returnderivativebound}, the sequence of maps $$L_j:=H_{\ta(a^{m_j})(y)} \circ \ta(a^{m_j-n_j}) \circ H_{\ta(a^{n_j})(y)}^{-1}\in \cG_1$$ has uniformly bounded derivatives and thus converges to a map $L\in  \cG_1$.
 Moreover, we have $(L_j)_*\Psi_{\ta(a^{n_j})(y)}$ is proportional to $\Psi_{\ta(a^{m_j})(y)} $   and thus
 $L_*\Psi_x$ is proportional to $\Psi_{x'}$.  Thus there exists a map $T= H_{x} \circ H_{x'}^{-1}\circ L
  \in \cG$ with $T(0) = H_x(x')$ and $T_*\Psi_x$ proportional to $\Psi_x$. 
Then arguing as before, as in Lemmas 6.1--6.3 in \cite{KS} following Section 5 of \cite{KS2}, we deduce that  $\Psi_x$ is supported on an orbit of a subgroup $H$ of $\cG$, and is the image of the Haar measure on $H$.  Thus $\mu^{\cW^{[\beta]}, {\cE_\beta}}_x$ is in the Lebesgue class of an embedded submanifold $N_x \subset \cW^{[\beta]}(x)$ for $\mu$-a.e. $x$.
This completes the proof for (2).
\end{proof}

\subsection{The proof of Theorem \ref{prop:dim2}} \label{subsec Proofpropdim2}
\hfill

Using the Lemma \ref{lem measurerigidity}, we have obtained the following dichotomy.
\begin{proposition}\label{prop:ergod} \label{prop:entorsmooth}
Let $\mu$ be a $P$-invariant, $P$-ergodic Borel probability measure on $M^\alpha$.  Consider a restricted root $\beta$ such that $\dim E^{[\beta]}_F(x) = 1$ for $\mu$-a.e. $x$. Then either:
\begin{enumerate}
\item  $\mu$ is  $G_{[\beta]}$-invariant, or
\item  $\Lg_{[\beta]}$ is one dimensional and for $\mu$-a.e. $x$, the leafwise measure
 $\mu^{\cW^{[\beta]},{\cE_\beta}}_x$ is in the Lebesgue class on a properly embedded $C^1$ curve $N_x\subset \cW^{[\beta]}(x)$ that is everywhere   transverse to the $G_{[\beta]}$-orbits.
\end{enumerate}
\end{proposition}

\begin{proof} 
To establish Proposition \ref{prop:entorsmooth}, we consider first the case that  $\mu^{G_{[\beta]},{\cE_\beta}}_x$ is not a single atom for $\mu$-a.e. $x$.
Then by Lemma \ref{lem measurerigidity}, we deduce that $\mu^{G_{[\beta]},{\cE_\beta}}_x$ is in the Lebesgue class of an embedded submanifold $N_x \subset G_{[\beta]}\cdot x$ for $\mu$-a.e. $x$.
From Proposition \ref{prop:moreentr}, it follows that  $\mu$ is  $G_{[\beta]}$-invariant.

We may thus assume from now on that  the leafwise measures $\mu^{G_{[\beta]},{\cE_\beta}}_x$ along $G_{[\beta]}$-orbits are Dirac masses.
By Lemma \ref{lem measurerigidity},  we conclude the leafwise measure $\mu^{\cW^{[\beta]},{\cE_\beta}}_x$ is in the Lebesgue class of an embedded submanifold $N_x \subset \cW^{[\beta]}( x)$ for almost every $x$.
 Note that the above discussions rely on the hypothesis $\dim E^{[\beta]}_F(x) = 1$ since we have used Lemma \ref{lem measurerigidity}.

By applying  Lemma \ref{lemma:entropy}  to $\gamma = \beta$, we deduce the following.
\begin{corollary}
For $\mu$-a.e.\ $x$,   the embedded submanifold $N_x$
has  codimension at most $1$ in $\cW^{[\beta]}(x)$.
\end{corollary}
\begin{proof}
In all of the following discussions, we will let $x$ denote a $\mu$-typical point.

By construction, the Lebesgue class of an embedded submanifold $N_x$ is equivalent to the leafwise measure $\mu^{\cW^{[\beta]},{\cE_\beta}}_x$. Then $\cN = \{ N_x \}_{x \in M}$ forms a measurable foliation invariant under the $A$-action. By definition, we see that $T_x N_x$ forms an $A$-invariant subspace of $E^{[\beta]}$.
By comparing the Lyapunov functionals of the Oseledects splitting of $T_x N_x$ and those of $E^{[\beta]}(x)$,  we see that there is a splitting $T_x N_x = E_{0}(x) \oplus E_{1}(x) \oplus E_{2}(x)$ where $E_{0}(x) = T_x N_x \cap E_F(x)$ and $E_{i}(x) \cap E_F(x) = \emptyset$ for $i \in \{1,2\}$. Moreover $E_i$ corresponds to Lyapunov functional $i \beta$.  If $\cN$ has codimension at least $2$, then by hypothesis $\dim E_F^{[\beta]} = 1$, we see that either $\dim E_1(x) < \dim \Lg_{\beta}$ or $\dim E_2(x) < \dim \Lg_{2\beta}$.

Let $\lambda^F$ denote the fiberwise Lyapunov exponent in $[\beta]$.
Since $\mu^{\cW^{[\beta]},{\cE_\beta}}_x$ is assumed smooth along a submanifold, an analogue of the Pesin entropy formula shows that
$	h_\mu(a\mid {\cE_\gamma} \vee  \cW^{[ \gamma ]})= 	\lambda(a)+ \gamma( a) ( \dim E_{1}(x) + 2	 \dim  E_{2}(x) )$
and
 $h_\mu(a\mid {\cE_\gamma }\vee \cF \vee \cW^{[ \gamma ]})= \lambda(a)$. 
Combining these contradicts the equality \eqref{eq:entgaplower}. Thus $\cN$ has codimension at most $1$ in $\cW^{[\beta]}$.
\end{proof}

We may deduce from the above corollary that $\dim \Lg_{[\beta]} = 1$, for otherwise
we must have $T_x(G_{[\beta]}\cdot x)\cap T_xN_x\neq\{0\}$ for $\mu$-a.e. $x$. This would imply that $\mu$ is $G_{[\beta]}$-invariant by Proposition \ref{prop:moreentr}.  By the same reason, we see that $N_x$  is transverse to  the $G_{[\beta]}$-orbits at a positive measure set of $x\in M^\alpha$.  By ergodicity and $A$-equivariance, this implies $N_x$ is transverse to the $G$-orbit at a full measure set of $x\in M^\alpha$.  Since transversality is an open condition, by acting by an element of $A$ that expands $\cW^{[\beta]}$-leaves, we conclude for $\mu$-a.e.\ $x$ that $N_x$ is  everywhere transverse to the $G_{[\beta]}$-orbits. This concludes the proof.
\end{proof}

We can now conclude the proof of Theorem \ref{prop:dim2}.

\begin{proof}[Proof of Theorem \ref{prop:dim2}]
It suffices to show that: if $\dim E^{[\beta]}_F < 2$, then $\dim \Lg_{[\beta]} = 1$.

If  $E^{[\beta]}_F = \{0\}$, then we would have $[\beta] \cap \Sigma \subset \Sigma_Q$ by Theorem \ref{T:BRHW}.
Now assume that  $\dim E^{[\beta]}_F = 1$. Then by Proposition \ref{prop:entorsmooth}, we must have $\dim \Lg_{[\beta]} = 1$.
This concludes the proof.
\end{proof}

\subsection{Proof of Lemma \ref{lemma:entropy}} \label{subsec finishinglemmaentropy}

In this subsection we give the proof of Lemma \ref{lemma:entropy}. Recall that $\mu, P, Q$ are as in Theorem \ref{prop:dim2}, and we let $\gamma \in \Sigma \setminus \Sigma_Q$ satisfy $\frac 1 2  \gamma  \notin \Sigma$. We emphasize that, unlike in Subsections \ref{subsec measurerigidityhomogenity} and \ref{subsec Proofpropdim2}, we make no assumptions on the dimension of $E^{[\gamma]}_F$.

We begin with the following version of Abramov-Rohlin theorem.
\begin{claim}
Given $a\in A$, we have
\aryst
h_\mu(a\mid {\cE_\gamma}) = h_\mu(a\mid {\cE_\gamma}\vee \cF) + h_{m_{G/\Gamma}} (a).
\earyst
\end{claim}
\begin{proof}

\def\sigA{\mathfrak A}
\def\sigB{\mathfrak B}
\def\sigC{\mathfrak C}

\def\sigD {\mathfrak D}
Let $\sigB$ denote the ambient (Borel) $\sigma$-algebra on $M^\alpha$.
Let $\sigC$ be the sub-$\sigma$-algebras of $\sigB$ generated by the partition $\cF$ into fibers and let $\sigD$ be the sub-$\sigma$-algebras of $\sigB$ generated by the $A'_\gamma$-ergodic decomposition  ${\cE_\gamma}$ of $\mu$.

Given sub-$\sigma$ algebras $\sigA,\sigC\subset \sigB$ with $\sigC$ $a$-invariant, write
\[
h_{ \mu}(a, \sigA\mid \sigC):=\sup _{\mathcal P\subset  \sigA} h _{ \mu}(a, \mathcal P|\sigC) \]
where the supremum is taken over all finite sub-$\sigma$-algebras (i.e. finite partitions)  $\mathcal P\subset  \sigA$.  With this notation, we have

\begin{enumerate}
	\item $h_\mu(a\mid {\cE_\gamma}) = h_{ \mu}(a,\sigB \mid \sigD)$,
	\item $ h_\mu(a\mid {\cE_\gamma}\vee \cF) = h_{ \mu}(a, \sigB\mid \sigD \vee \sigC)$, and
	\item  $h_{m_{G/\Gamma}}(a) = h_{ \mu}(a, \sigC)$.
\end{enumerate}

Recall that $\pi\colon M^\alpha\to G/\Gamma$ is the canonical projection.  Then $\sigC$ is the $\sigma$-algebra of level sets under $\pi$.
Recall that for every non-identity element $a\in A$, the $a$-action on $G/\Gamma$ is ergodic with respect to the Haar measure $m_{G/\Gamma}$.
Given an $A'_\gamma$-ergodic component  $\mu^\sigD_x$ of $\mu$, we have $\pi_*\mu^\sigD_x= m_{G/\Gamma}$ for $\mu$-a.e.\ $x$.
In particular,  given $C\in \sigC$ we have
 \aryst
 \mu^\sigD_x(\pi^{-1} (\pi(C))) = \pi_*(\mu^\sigD_x)(\pi(C))= m_{G/\Gamma}(\pi(C)).
 \earyst
Thus, for $a\in A'_\gamma$ and any finite partition $\mathcal P$ of $G/\Gamma$,
\begin{align*}
h_\mu(a, \pi^{-1}(\mathcal P)\mid \sigD)
&= \int -  \sum_{C\in \mathcal P}
 \mu^\sigD_x(\pi^{-1} C) \log \mu^\sigD_x(\pi^{-1} C) \,d\mu(x) \\
&= \int - \sum_{C\in \mathcal P}
 m_{G/\Gamma}(C) \log m_{G/\Gamma}(C)\,
d\mu(x) \\
&= h_{m_{G/\Gamma}}(a, \mathcal P)
\end{align*}
Then, by the definition of $\sigC$, we have
\aryst
h_{ \mu}(a, \sigC)= h_{m_{G/\Gamma}}( a)= h_{ \mu}(a, \sigC\mid \sigD).
\earyst
Also, since $\sigC\subset \sigB$,
\aryst
h_{\mu}(a,\sigB \mid \sigD)=h_{\mu}(a,\sigB \vee \sigC\mid \sigD).
\earyst
Finally, we have (see e.g.\cite[Fact 4.1.7]{Do})
\aryst
h_{ \mu}(a,\sigB\vee\sigC \mid \sigD) = h_{ \mu}(a,\sigB\mid \sigC\vee \sigD) + h_{ \mu}(a,\sigC\mid \sigD).
\earyst
The result then follows from the above equalities.
\end{proof}

\begin{proof}[Proof of Lemma \ref{lemma:entropy}]

All three terms in \eqref{eq:entgaplower} are positively proportional and restrict to linear function on the half-space  $\gamma(a)>0$.   In particular, \eqref{eq:entgaplower} holds for $a$ if it holds for at least one $a'\in A$ with  $\gamma(a')>0$.
It is therefore enough to verify \eqref{eq:entgaplower}  for any
$a\in A$ such  that
\begin{enumerate}
	\item $\gamma(a)>0$ for all $\gamma \in \Sigma \setminus \Sigma_Q$, and
	\item  $\gamma(a)\le 0$ for all $\gamma \in \Sigma_Q$.
\end{enumerate}

Partition the coarse Lyapunov functionals $\widehat \cL(A, \mu)= \widehat \cL_1\sqcup\widehat \cL_2$ so that $\chi\in \widehat \cL_1$  if and only if  $\chi =[\beta']$ for some $\beta' \in \Sigma \setminus \Sigma_Q$.
Given $\chi\in \widehat \cL_2$, we have either $\chi(a)\le 0$ in which case $\cW^\chi$ contributes no entropy, or  $\cW^{\chi}_{F}(x) = \cW^\chi(x)$ for almost every $x$ and so $$h_ \mu(a\mid {\cE_\gamma}\vee \cW^\chi)= h_ \mu(a\mid {\cE_\gamma}\vee \cF\vee \cW^\chi). $$

		We have
\begin{align*}
	h_\mu(a\mid &{\cE_\gamma})
	= h_\mu(a\mid {\cE_\gamma}\vee \cF) + h_{{m_{G/\Gamma}}} (a)\\
	&=\sum_{\substack{\chi\in \widehat \cL(A, \mu) \\ \chi(a)>0}} h_\mu(a\mid {\cE_\gamma}\vee \cF\vee \cW^\chi)+
 \sum_{\substack{ [\gamma]: \gamma \in  \Sigma \setminus \Sigma_Q }} h_ {m_{G/\Gamma}}(a\mid W^{[\gamma]}_G ). 
\end{align*}
On the other hand,
\begin{align*}
	h_\mu(a\mid &{\cE_\gamma})
	=\sum_{\substack{\chi\in \widehat \cL(A, \mu) \\ \chi(a)>0}} h_ \mu(a\mid {\cE_\gamma}\vee \cW^\chi). 
\end{align*}
Thus
\begin{align*}
	\sum_{\substack{\chi\in \widehat \cL_1}} h_ \mu(a\mid {\cE_\gamma}\vee \cW^\chi) &+
		\sum_{\substack{\chi\in \widehat \cL_2\\ \chi(a)>0}} h_ \mu(a\mid {\cE_\gamma}\vee \cW^\chi)\\
		&\ =\sum_{\substack{\chi\in \widehat \cL_1}} h_\mu(a\mid {\cE_\gamma}\vee \cF\vee \cW^\chi)+
		\sum_{\substack{\chi\in \widehat \cL_2\\ \chi(a)>0}} h_\mu(a\mid {\cE_\gamma}\vee \cF\vee \cW^\chi)\\&\quad +
 \sum_{\substack{ [\gamma] : \gamma \in  \Sigma \setminus \Sigma_Q }} h_ {m_{G/\Gamma}}(a\mid W^{[\gamma]}_G  ).
\end{align*}
Canceling common $\widehat \cL_2$ terms,
\begin{align*}
\sum_{\substack{ [\gamma] : \gamma \in  \Sigma \setminus \Sigma_Q }} \bigl(h_ \mu(a\mid {\cE_\gamma}\vee \cW^{[\gamma]})
	-h_\mu(a\mid {\cE_\gamma}\vee \cF\vee \cW^{[\gamma]})\bigr)=
	\sum_{\substack{ [\gamma] : \gamma \in  \Sigma \setminus \Sigma_Q }} h_ {m_{G/\Gamma}}(a\mid W^{[\gamma]}_G  ).
\end{align*} 
 Combined with \eqref{basicLYbound} and the fact that
\aryst
 h_ {m_{G/\Gamma}}(a\mid W^{[\gamma]}_G  ) =  	 \gamma( a) ( \dim \Lg_\gamma + 2	 \dim \Lg_{2 \gamma } ),
\earyst
 the result follows.
\end{proof}

\subsection{A new dimension bound and applications} \label{subsec anewdimbound}

For $[\beta]\in\widehat{\Sigma}$, we  let $r_0([\beta]) = 1$ if $\dim \Lg_{[\beta]} = 1$ and let $r_0([\beta]) = 2$ otherwise.
If $Q\subset G$ is a standard parabolic subgroup, we let
$$r_0(Q)=\sum_{[\beta]\in\widehat{\Sigma}\sm\widehat{\Sigma}_Q} r_0([\beta]).$$

\begin{prop} \label{prop dim geq r0}
Let $P$, $\mu$ and $Q$ be as in Theorem \ref{prop:dim2}. Then
$$\dim(M)\ge r_0(Q).$$
\end{prop}

\begin{proof}
By Theorem \ref{prop:dim2}, we have
$$\dim(M) = \sum_{\chi \in \widehat \cL(A, \mu)} \dim E^{\chi}_F \geq \sum_{[\beta]\in\widehat{\Sigma}\sm\widehat{\Sigma}_Q} \dim E^{[\beta]}_F
\ge\sum_{[\beta]\in\widehat{\Sigma}\sm\widehat{\Sigma}_Q} r_0([\beta])=r_0(Q).$$
\end{proof}

We also let
$$r_0(G)=\min_{Q\subsetneq G \mbox{\tiny{ standard parabolic}}} r_0(Q).$$
The following corollary of Proposition \ref{prop dim geq r0} is a direct consequence of the definition of $s(G)$.

\begin{cor} \label{cor s geq r0}
We have $s(G) \geq r_0(G)$.
\end{cor}

We are now prepared to prove Theorem \ref{thm valuesofs}(1) and Theorem \ref{T:main-2}.

\begin{proof}[Proof of Theorem \ref{thm valuesofs}(1)]
Let $G$ be either a connected complex semisimple Lie group without simple factors of rank 1, or the group EII. Then $\dim\Lg_{\chi}\le2$ for every $\chi\in\widehat\Sigma$ (for complex groups we always have $\dim\Lg_{\chi}=2$; for the group EII, see \cite[page 534]{He}). This implies that
$r_0([\beta]) = \dim \Lg_{[\beta]}$ for every $[\beta]\in\widehat{\Sigma}$, and hence $r_0(Q)=\dim G/Q$ for every  standard parabolic
subgroup. It follows that $r_0(G)=v(G)$, and the proof is concluded by applying Corollary \ref{cor s geq r0} and Lemma \ref{L:s-ineq}.
\end{proof}

\begin{remark}
The above proof works equally well for any group such that $\dim\Lg_{\chi}\le2$ for every $\chi\in\widehat\Sigma$.
Such groups also include $\SO^+(n,n+2)$ and $\SU(n,n)$. However, the results for the latter two groups have been incorporated into \eqref{E:s-SO} and \eqref{E:s-SU}, whose proofs will be given by combining the superrigidity approach.
\end{remark}

\begin{proof}[Proof of Theorem \ref{T:main-2}]
Let $G$, $\Gamma$, $M$ and $\alpha$ be as in Theorem \ref{T:main-2}. As noted above, we have $\dim\Lg_{\chi}\le2$ for every $\chi\in\widehat\Sigma$.
In view of Theorem \ref{T:general}, it suffices to prove that every $P$-invariant Borel probability measure on $M^\alpha$ is $G$-invariant. Assume to the contrary that there exists a $P$-invariant Borel probability measure $\mu$ on $M^\alpha$ which is not $G$-invariant. Without loss of generality, we may assume that $\mu$ is $P$-ergodic. Let $Q=\Stab_G(\mu)$, which is a proper standard parabolic subgroup of $G$. By Theorem \ref{prop:dim2},
we have $\widehat\Sigma\sm\widehat\Sigma_Q\subset\widehat\cL_{F}(A,\mu)$ and $\dim E^{\chi}_F(x) \geq \dim \Lg_{\chi}$ for every $\chi\in\widehat\Sigma\sm\widehat\Sigma_Q$. It follows that
\begin{align*}
v(G)&=\dim(M)=\sum_{\chi \in \widehat \cL_F(A, \mu)} \dim E^{\chi}_F\ge\sum_{\chi \in\widehat\Sigma\sm\widehat\Sigma_Q} \dim E^{\chi}_F\\
&\ge\sum_{\chi \in\widehat\Sigma\sm\widehat\Sigma_Q}\dim \Lg_{\chi}=\dim G/Q\ge v(G).
\end{align*}
So we must have $\widehat\cL_{F}(A,\mu)=\widehat\Sigma\sm\widehat\Sigma_Q$, namely, every $\lambda\in\cL_F(A,\mu)$ is a positive multiple of some root in $\Sigma\sm\Sigma_Q$. However, since the action $\alpha$ preserves a smooth volume density, we have
$$\sum_{\lambda\in\cL_F(A,\mu)}(\dim E_{F, \lambda})\lambda=0.$$
This leads to a contradiction.
\end{proof}

As we have mentioned at the end of Subsection \ref{subsec twoapproaches}, the measure rigidity approach developed in this section can be combined with the superrigidity approach in Section \ref{sec: superrigidity}
to get better estimates. In the rest of this section, we will illustrate this point by proving  \eqref{E:s-SL} and \eqref{E:s-SO}. We will see more applications of this combined approach in Section \ref{sec: comination}.

\begin{proof}[Proof of \eqref{E:s-SL}]
Let $G=\SL(n,\HH)$, $n\ge3$. We can take $\La$ to be the space of real diagonal matrices of trace zero. Then
$$\Sigma=\{e_i-e_j:1\le i\ne j\le n\},$$
where $e_i\in\La^*$ is the functional sending a diagonal matrix to its $i$-th diagonal entry.
We have that
$$\dim\Lg_{e_i-e_j}=4.$$
Thus
$$r_0([e_i-e_j])=2.$$
The set of simple roots $\Pi=\{\alpha_1,\ldots,\alpha_{n-1}\}$ can be chosen as
$$\alpha_1=e_1-e_2, \quad \ldots, \quad \alpha_{n-1}=e_{n-1}-e_n.$$
We need to prove that if $\alpha$ is an action of a lattice $\Gamma\subset G$ on a compact manifold $M$ such that there is a $P$-invariant Borel probability measure on $M^\alpha$ that is not $G$-invariant, then $\dim(M)$ is not less than the right side hand of \eqref{E:s-SL}.
Without loss of generality, we may assume that $\mu$ is $P$-ergodic. Let $Q=\Stab_G(\mu)$. Then $\Pi_Q\ne\Pi$. There are 2 cases.

\textbf{Case 1.} Suppose that $\Pi_Q$ has a connected component $\Delta$ with $|\Delta|\ge2$, say $\Delta=\{\alpha_p,\alpha_{p+1},\ldots,\alpha_q\}$, where $1\le p<q\le n-1$ and $(p,q)\ne(1,n-1)$. Then $\Lg_\Delta\cong\Lsl(q-p+2,\HH)$. By Proposition \ref{P:real-1}, we have
$$\dim(M)\ge |\Pi\sm\Delta|n(\Lg_\Delta)=4(n-q+p-2)(q-p+2)\ge 4n-4,$$
which is not less than the right side hand of \eqref{E:s-SL}.

\textbf{Case 2.} Suppose every connected component of $\Pi_Q$ consists of a single root. Then $\Sigma\sm\Sigma_Q=-(\Sigma^+\sm\Pi_Q)$ and $|\Pi_Q|\le[n/2]$. It follows from Proposition \ref{prop dim geq r0} that
$$\dim(M)\ge r_0(Q)=2|\Sigma\sm\Sigma_Q|=2|\Sigma^+\sm\Pi_Q|\ge n(n-1)-2[n/2],$$
which is again not less than the right side hand of \eqref{E:s-SL}.
This proves \eqref{E:s-SL}.
\end{proof}

\begin{proof}[Proof of \eqref{E:s-SO}]
Let $G=\SO^+(m,n)$ with $m-2\ge n\ge2$. The root system $\Sigma$ is of type $B_n$, namely, there exists an orthogonal basis $\{e_1,\ldots,e_n\}$ of $\La^*$ with $\|e_1\|=\cdots=\|e_n\|$ such that
$$\Sigma=\{\pm e_i\pm e_j:1\le i\ne j\le n\}\cup\{\pm e_i:1\le i\le n\}.$$
Moreover, the dimensions of the restricted root spaces are
$$\dim\Lg_{\pm e_i\pm e_j}=1, \qquad \dim\Lg_{\pm e_i}=m-n.$$
Thus
$$r_0([\pm e_i\pm e_j])=1, \qquad r_0([\pm e_i])=2.$$
The set of simple roots $\Pi=\{\alpha_1,\ldots,\alpha_n\}$ can be chosen as
$$\alpha_1=e_1-e_2, \quad \ldots, \quad \alpha_{n-1}=e_{n-1}-e_n,  \quad \alpha_n=e_n.$$
It is straightforward to verify that \eqref{E:s-SO} is equivalent to
$$s(G)\ge\min\{m+n-2, n(n+3)/2, (n-1)(n+6)/2\}.$$
To prove the later inequality, we need to show that if $\alpha$ is an action of a lattice $\Gamma\subset G$ on a compact manifold $M$ such that there is a $P$-invariant Borel probability measure on $M^\alpha$ that is not $G$-invariant, then
$$\dim(M)\ge\min\{m+n-2,n(n+3)/2, (n-1)(n+6)/2\}.$$
Without loss of generality, we may assume that $\mu$ is $P$-ergodic.
Let $Q=\Stab_G(\mu)$. Then $\Pi_Q\ne\Pi$. Let $k\in\{1,\ldots,n\}$ be the largest index such that $\alpha_k\notin\Pi_Q$. There are 2 cases.

\textbf{Case 1.} Suppose $1\le k\le n-2$. Then $\Delta=\{\alpha_{k+1},\ldots,\alpha_n\}$ is a connected component of $\Pi_Q$ with $|\Delta|=n-k\ge2$,
and we have $\Lg_\Delta\cong\Lso(m-k,n-k)$ (a real simple Lie algebra is determined up to isomorphism by its restricted root system and the dimensions of the restricted root spaces, see,
for example, \cite[Chapter X, Exercise F.9(b)]{He}). By Proposition \ref{P:real-1}, we have
$$\dim(M)\ge kn(\Lg_\Delta)=k(m+n-2k)\ge m+n-2.$$

\textbf{Case 2.} Suppose $k\ge n-1$. It follows from $\alpha_k\notin\Pi_Q$ that $\Sigma_Q$ does not contain the following $k(4n-3k+1)/2$ roots:
$$\{-e_i:1\le i\le k\}\cup\{-e_i-e_j:1\le i<j\le k\}\cup\{-e_i\pm e_j:1\le i\le k, k+1\le j\le n\}.$$
By Proposition \ref{prop dim geq r0}, we have
$$\dim(M)\ge r_0(Q)\ge k(4n-3k+1)/2+k\ge\min\{n(n+3)/2, (n-1)(n+6)/2\}.$$
This completes the proof.
\end{proof}

\section{Combination of the two approaches: $\SU(m,n)$, $\Sp(m,n)$ and $\SO^*(2n)$.} \label{sec: comination}

In this section, we prove \eqref{E:s-SU}--\eqref{E:s-SO*}.
To unify the notation, we let
$$\epsilon_G=\begin{cases}
  2, & \text{if $G=\SU(m,n)$ with $m\ge n\ge2$;} \\
  4, & \text{if $G=\Sp(m,n)$ with $m\ge n\ge2$, or $\SO^*(2n)$ with $n\ge4$.}
\end{cases}$$
Moreover, it will be convenient to rewrite \eqref{E:s-SU} as the following two inequalities:
\begin{align}
s(\SU(n,n)) & \ge\min\{4n-3,n^2\}, && n\ge2, \label{E:s-SU1} \\
s(\SU(m,n)) & \ge\min\{2m+2n-3,n(n+1)\}, && m>n\ge2, \label{E:s-SU2}
\end{align}
and rewrite \eqref{E:s-SO*} as the following two inequalities:
\begin{align}
s(\SO^*(4n)) & \ge\min\{8n-7,n^2\}, && n\ge2, \label{E:s-SO*1} \\
s(\SO^*(4n+2)) & \ge\min\{8n-3,n(n+1)\}, && n\ge2. \label{E:s-SO*2}
\end{align}
Then \eqref{E:s-SU1} and \eqref{E:s-SO*1} can be unified as
\begin{equation}\label{E:s-uni1}
s(G)\ge\min\{2\epsilon_G(n-1)+1,n^2\}, \qquad G=\SU(n,n) \text{ or } \SO^*(4n) \text{ with } n\ge2,
\end{equation}
and \eqref{E:s-SU2}, \eqref{E:s-Sp} and \eqref{E:s-SO*2} can be unified as
\begin{equation}\label{E:s-uni2}
s(G)\ge\min\{\epsilon_G(m+n-2)+1,n(n+1)\}, \qquad G=\begin{cases}
\SU(m,n) \ (m>n\ge2), \text{ or}\\
\Sp(m,n) \ (m\ge n\ge2), \text{ or}\\
\SO^*(2m+2n) \ (m=n+1).
\end{cases}
\end{equation}

For the groups $\SU(m,n)$, $\Sp(m,n)$, $\SO^*(4n)$ and $\SO^*(4n+2)  \ (m\ge n\ge2)$, the restricted root system $\Sigma$ has type $C_n$ or $BC_n$. More precisely, there is an orthogonal basis $\{e_1,\ldots,e_n\}$ of $\La^*$ with $\|e_1\|=\cdots=\|e_n\|$ such that:
\begin{itemize}
  \item If $G=\SU(n,n)$, $\Sp(n,n)$ or $\SO^*(4n)$, where $n\ge2$, then
$$\Sigma=\{\pm e_i\pm e_j:1\le i\ne j\le n\}\cup\{\pm 2e_i:1\le i\le n\}.$$
We have
$$\dim\Lg_{\pm e_i\pm e_j}=\epsilon_G, \qquad \dim\Lg_{\pm 2e_i}=\begin{cases}
  1, & G=\SU(n,n) \ \text{ or } \ \SO^*(4n);\\
  3, & G=\Sp(n,n).
\end{cases}$$
In this case, we choose $\Pi=\{\alpha_1,\ldots,\alpha_n\}$ as
$$\alpha_1=e_1-e_2, \quad \ldots, \quad \alpha_{n-1}=e_{n-1}-e_n,  \quad \alpha_n=2e_n.$$
  \item If $G=\SU(m,n) \ (m>n\ge2)$, $\Sp(m,n) \ (m>n\ge2)$ or $\SO^*(4n+2) \ (n\ge2)$, then
$$\Sigma=\{\pm e_i\pm e_j:1\le i\ne j\le n\}\cup\{\pm e_i,\pm 2e_i:1\le i\le n\}.$$
We have
$$\dim\Lg_{\pm e_i\pm e_j}=\epsilon_G, \qquad \dim\Lg_{\pm e_i}=\epsilon_G(m-n), \qquad \dim\Lg_{\pm 2e_i}=\begin{cases}
  1, & G=\SU(m,n) \ \text{ or } \ \SO^*(4n+2);\\
  3, & G=\Sp(m,n).
\end{cases}$$
In this case, we choose $\Pi=\{\alpha_1,\ldots,\alpha_n\}$ as
$$\alpha_1=e_1-e_2, \quad \ldots, \quad \alpha_{n-1}=e_{n-1}-e_n,  \quad \alpha_n=e_n.$$
\end{itemize}
Note that in both cases, for $\beta\in\Sigma$ we have
\begin{equation}\label{E:r0}
r_0([\beta])=\begin{cases}
  1, & \text{if } \ G\in\{\SU(n,n), \SO^*(4n)\} \ \text{ and } \ \beta=\pm 2e_i;\\
  2, & \text{otherwise}.
\end{cases}
\end{equation}
To prove \eqref{E:s-uni1} and \eqref{E:s-uni2}, we need the following result, which strengthens Proposition \ref{P:real-1} in special cases.

\begin{proposition}\label{P:real-2}
Suppose $G$ is one of the following groups:
\begin{itemize}
  \item $\SU(m,n)$ with $m\ge n\ge3$ and $(m,n)\ne(3,3)$;
  \item $\Sp(m,n)$ with $m\ge n\ge3$;
  \item $\SO^*(2m+2n)$ with $n\ge3$, $m\in\{n,n+1\}$ and $(m,n)\ne(3,3)$.
\end{itemize}
Let $P\subset G$ be the minimal standard parabolic subgroup,
$Q\subset G$ be the standard parabolic subgroup with $\Pi_Q=\Delta:=\{\alpha_2,\ldots,\alpha_n\}$,
$\Gamma\subset G$ be a lattice, $M$ be a compact manifold, $\alpha:\Gamma\to\Diff^{\reg}(M)$ be an action, and $\mu$ be a $P$-invariant $P$-ergodic Borel probability measure on $M^\alpha$ such that $\Stab_G(\mu)=Q$. Then
  $$\dim(M)\ge n(\Lg_\Delta)+1,$$
where $\Lg_\Delta$ and $n(\Lg_\Delta)$ are as in Proposition \ref{P:real-1}.
\end{proposition}

We first prove two lemmas.

\begin{lemma}\label{L:root-space}
Let $\Lg$ be a real semisimple Lie algebra, $\La\subset\Lg$ be a maximal $\RR$-diagonalizable subalgebra, and $\Sigma=\Sigma(\Lg,\La)$. Let $\beta\in\Sigma$ and $\gamma\in\Sigma\cup\{0\}$ be such that $\langle\beta,\beta+2\gamma\ra\ge0$. Then for every $\bx\in\Lg_{\beta}\sm\{0\}$, we have $[\bx,\Lg_\gamma]=\Lg_{\beta+\gamma}$.
\end{lemma}

\begin{proof}
Let $\by\in\Lg_{-\beta}\sm\{0\}$ and $\bh\in\La\sm\{0\}$ be such that $\Lsl_\bx:=\spann\{\bh,\bx,\by\}$ is a subalgebra isomorphic to $\Lsl(2,\RR)$, and such that the linear map determined by
$$\bh\mt\begin{pmatrix}
  1 & 0 \\ 0 & -1
\end{pmatrix}, \qquad \bx\mt\begin{pmatrix}
  0 & 1 \\ 0 & 0
\end{pmatrix}, \qquad \by\mt\begin{pmatrix}
  0 & 0 \\ 1 & 0
\end{pmatrix}$$
is an isomorphism (see \cite[Proposition 6.52]{Kn}). Consider the adjoint representation of $\Lsl_\bx$ on $\bigoplus_{k\in\ZZ}\Lg_{k\beta+\gamma}$.
Note that $\beta(\bh)=2$, $\gamma(\bh)=\frac{2\langle\beta,\gamma\ra}{\langle\beta,\beta\ra}$. So $\Lg_{k\beta+\gamma}$ (if it is nonzero) is the eigenspace for $\ad(\bh)$ with eigenvalue $2k+\frac{2\langle\beta,\gamma\ra}{\langle\beta,\beta\ra}$. In particular, $\Lg_\gamma$ is an eigenspace with eigenvalue not less than $-1$. It follows that the map $\ad(\bx):\Lg_\gamma\to\Lg_{\beta+\gamma}$ is surjective (see \cite[Ch. 8, Sec. 1.2, Corollary to Proposition 2]{Bou2}). This means $[\bx,\Lg_\gamma]=\Lg_{\beta+\gamma}$.
\end{proof}

As in \cite[Section I.8]{Kn}, the Lie algebras
\begin{equation}\label{E:rep}
\begin{cases}
\Lsu(m,n), & m\ge n\ge2, \ (m,n)\ne(2,2);\\
\Lsp(m,n), & m\ge n\ge2; \\
\Lso^*(2n), & n\ge5
\end{cases}
\end{equation}
are defined to be subalgebras of $\Lgl(m+n,\CC)$, $\Lgl(m+n,\HH)$ and $\Lgl(2n,\CC)$, respectively. If $\Lg$ is one of them, its definition gives an $n(\Lg)$-dimensional real representation by composing with the natural inclusion $\Lgl(d,\CC)\hookrightarrow\Lgl(2d,\RR)$ or $\Lgl(d,\HH)\hookrightarrow\Lgl(4d,\RR)$. Let us call it the standard representation and denote it by $\st:\Lg\to\Lgl(n(\Lg),\RR)$. The next lemma is as follows.

\begin{lemma}\label{L:rep}
Suppose $\Lg$ is one of the Lie algebras in \eqref{E:rep}. Then every $n(\Lg)$-dimensional nontrivial real representation of $\Lg$ is equivalent to the standard representation.
\end{lemma}

\begin{proof}
We first review some facts on real representations of a real semisimple Lie algebra $\Lg$ (see \cite[\S 8]{On}). Fix a set of fundamental weights $\{\varpi_1,\ldots,\varpi_\ell\}$ for the complexification $\Lg_\CC$ of $\Lg$, and let
$\Lambda=\{\sum_{i=1}^\ell k_i\varpi_i:k_i\in\ZZ_{\ge0}\}$
be the semigroup of dominant integral weights. With respect to the real form $\Lg$ of $\Lg_\CC$, the following two invariants can be introduced:
\begin{itemize}
  \item An involution $s_0:\Lambda\to\Lambda$ corresponding to an involution of the Dynkin diagram of $\Lg_\CC$. Denote $\Lambda^{s_0}=\{\lambda\in \Lambda:s_0(\lambda)=\lambda\}$.
  \item A homomorphism $\epsilon:\Lambda^{s_0}\to\{\pm1\}$, called the Cartan index. Denote $\ker(\epsilon)=\{\lambda\in\Lambda^{s_0}:\epsilon(\lambda)=1\}$.
\end{itemize}
For every $\lambda\in \Lambda$, fix an irreducible complex representation $\rho_\lambda:\Lg_\CC\to\Lgl(d_\lambda,\CC)$ with highest weight $\lambda$. Then, for an irreducible real representation $\phi:\Lg\to\Lgl(d,\RR)$, one has:
\begin{itemize}
  \item[(i)] If the complexification $\phi_\CC:\Lg_\CC\to\Lgl(d,\CC)$ of $\phi$ is irreducible, then $\phi_\CC\cong\rho_\lambda$ for some $\lambda\in\ker(\epsilon)$;
  \item[(ii)] If $\phi_\CC$ is reducible, then $\phi\cong\rho_\lambda^\RR$ for some $\lambda\in\Lambda\sm\ker(\epsilon)$, where $\rho_\lambda^\RR$ is the composition $\Lg\stackrel{\rho_\lambda|_\Lg}{\longrightarrow}\Lgl(d_\lambda,\CC)\hookrightarrow\Lgl(2d_\lambda,\RR)$. Moreover, $\rho_{s_0(\lambda)}^\RR\cong\rho_\lambda^\RR$.
\end{itemize}

Now let $\Lg$ be one of the Lie algebras in \eqref{E:rep}, and let $\phi_0:\Lg\to\Lgl(n(\Lg),\RR)$ be a nontrivial representation. Note that $\phi_0$ and $\st$ are irreducible. We show that $\phi_0\cong\st$ by proving the following two statements:
\begin{itemize}
  \item If $\phi:\Lg\to\Lgl(d,\RR)$ is a nontrivial irreducible representation such that $\phi_\CC$ is irreducible, then $d>n(\Lg)$;
  \item There exists a unique (up to equivalence) irreducible representation $\phi:\Lg\to\Lgl(n(\Lg),\RR)$ such that $\phi_\CC$ is reducible.
\end{itemize}
In view of (i) and (ii) above, it suffices to prove:
\begin{itemize}
  \item[(1)] $d_\lambda>n(\Lg)$ for every $\lambda\in\ker(\epsilon)\sm\{0\}$;
  \item[(2)] $\min\{d_\lambda:\lambda\in\Lambda\sm\{0\}\}=n(\Lg)/2$, and $d_\lambda=d_{\lambda'}=n(\Lg)/2$ implies $\lambda'\in\{\lambda,s_0(\lambda)\}$.
\end{itemize}
Note that $\Lg_\CC$ is of type $A_\ell$, $C_\ell$ or $D_\ell$. In each case, we assume that the fundamental weights $\varpi_1,\ldots,\varpi_\ell$ are enumerated as in \cite[Appendix C.1]{Kn}.

\textbf{Case 1.} Suppose $\Lg=\Lsu(m,n)$ with $m\ge n\ge2$ and $(m,n)\ne(2,2)$. Then $\Lg_\CC\cong\Lsl(\ell+1,\CC)$ is of type $A_\ell$, where $\ell=m+n-1\ge4$.
By \cite[Example 7.1 and Proposition 4.4]{On}, the involution $s_0$ of $\Lambda$ corresponds to the nontrivial automorphism of the Dynkin diagram of type $A_\ell$, namely $s_0(\varpi_i)=\varpi_{\ell+1-i}$ for $1\le i\le \ell$. So
$$\Lambda^{s_0}=\left\{\sum_{i=1}^\ell k_i\varpi_i:k_i\in\ZZ_{\ge0}, k_i=k_{\ell+1-i}\right\}.$$
A straightforward computation using the Weyl dimension formula shows that $d_{\varpi_i}=\binom{\ell+1}{i}$ and $d_{\varpi_1+\varpi_\ell}=\ell(\ell+2)$ (see \cite[Proposition 13.2]{Ca} and \cite[\S1, Example 24]{SK}). By \cite[\S1, Corollary 13]{SK}, for $\lambda\in\Lambda^{s_0}\sm\{0\}$ we have
$$d_\lambda\ge\min(\{d_{\varpi_1+\varpi_\ell}, d_{\varpi_2}+1,d_{\varpi_{\ell-1}}+1\}\cup\{d_{\varpi_i}:3\le i\le \ell-2\})>2\ell+2=n(\Lg).$$
So (1) holds. On the other hand, $d_{\varpi_1}=d_{\varpi_\ell}=\ell+1=n(\Lg)/2$, and for $\lambda\in\Lambda\sm\{0,\varpi_1,\varpi_\ell\}$ we have
$$d_\lambda\ge\min\{d_{\varpi_1}+1,d_{\varpi_2},\ldots,d_{\varpi_{\ell-1}},d_{\varpi_\ell}+1\}>\ell+1=n(\Lg)/2.$$
Since $s_0(\varpi_1)=\varpi_\ell$, (2) also holds.

\textbf{Case 2.} Suppose $\Lg=\Lsp(m,n)$, $m\ge n\ge2$. Then $\Lg_\CC\cong\Lsp(\ell,\CC)$ is of type $C_\ell$, where $\ell=m+n\ge4$. Since the Dynkin diagram of type $C_\ell$ has only trivial automorphism, the involution $s_0$ is trivial. So $\Lambda^{s_0}=\Lambda$. By \cite[Example 7.4]{On}, for $\lambda=\sum_{i=1}^\ell k_i\varpi_i\in\Lambda$, the Cartan index is given by
$$\epsilon(\lambda)=(-1)^{\sum_{i=1}^{[(\ell+1)/2]}k_{2i-1}}.$$
The Weyl dimension formula gives $d_{\varpi_i}=\frac{\ell-i+1}{\ell+1}\binom{2\ell+2}{i}$ and $d_{2\varpi_1}=\ell(2\ell+1)$ (see \cite[Proposition 13.23]{Ca} and \cite[\S1, Example 26]{SK}). It follows that for $\lambda\in\Lambda\sm\{0,\varpi_1\}$, we have
$$d_\lambda\ge\min\{d_{2\varpi_1},d_{\varpi_2},\ldots,d_{\varpi_\ell}\}>4\ell=n(\Lg).$$
Since $\varpi_1\notin\ker(\epsilon)$ and $d_{\varpi_1}=2\ell=n(\Lg)/2$, the above inequality implies (1) and (2).

\textbf{Case 3.} Suppose $\Lg=\Lso^*(2n)$ with $n\ge5$. Then $\Lg_\CC\cong\Lso(2n,\CC)$ is of type $D_n$. By \cite[Example 7.2 and Proposition 4.4]{On}, we have
\begin{itemize}
  \item If $n$ is even, the involution $s_0$ is trivial, and hence $\Lambda^{s_0}=\Lambda$;
  \item If $n$ is odd, $s_0$ corresponds to the nontrivial automorphism of the Dynkin diagram of type $D_n$, namely $s_0(\varpi_i)=\varpi_i$ for $1\le i\le n-2$, $s_0(\varpi_{n-1})=\varpi_n$, and $s_0(\varpi_n)=\varpi_{n-1}$. So
$$\Lambda^{s_0}=\left\{\sum_{i=1}^n k_i\varpi_i:k_i\in\ZZ_{\ge0}, k_{n-1}=k_n\right\}.$$
\end{itemize}
In both cases, for $\lambda=\sum_{i=1}^n k_i\varpi_i\in\Lambda^{s_0}$, the Cartan index is given by
$$\epsilon(\lambda)=(-1)^{\sum_{i=1}^{[n/2]}k_{2i-1}}.$$
We have $d_{\varpi_i}=\binom{2n}{i} \ (1\le i\le n-2)$, $d_{\varpi_{n-1}}=d_{\varpi_n}=2^{n-1}$, $d_{2\varpi_1}=(n+1)(2n-1)$ and $d_{\varpi_{n-1}+\varpi_n}=\binom{2n}{n-1}$ (see \cite[\S1, Example 28]{SK} and \cite[Problem V.9]{Kn}). It follows that if $n\ge6$, then for $\lambda\in\Lambda\sm\{0,\varpi_1\}$, we have
$$d_\lambda\ge\min\{d_{2\varpi_1},d_{\varpi_2},\ldots,d_{\varpi_n}\}>4n=n(\Lg).$$
Since $\varpi_1\notin\ker(\epsilon)$, this implies that (1) holds for $n\ge6$. If $n=5$, then $$\ker(\epsilon)=\{k_1\varpi_1+k_2\varpi_2+k_3\varpi_3+k_4(\varpi_4+\varpi_5):k_i\in\ZZ_{\ge0}, k_1+k_3 \text{ is even}\}.$$
It follows that for $\lambda\in\ker(\epsilon)\sm\{0\}$, we have
$$d_\lambda\ge\min\{d_{2\varpi_1},d_{\varpi_2},d_{\varpi_3},d_{\varpi_4+\varpi_5}\}>20=n(\Lg).$$
So (1) also holds for $n=5$. On the other hand, for any $n\ge5$, we have $d_{\varpi_1}=2n=n(\Lg)/2$ and
$$d_\lambda\ge\min\{d_{2\varpi_1},d_{\varpi_2},\ldots,d_{\varpi_n}\}>2n=n(\Lg)/2$$
for every $\lambda\in\Lambda\sm\{0,\varpi_1\}$. So (2) also holds.
\end{proof}

\begin{proof}[Proof of Proposition \ref{P:real-2}]
Assume to the contrary that $\dim(M)\le n(\Lg_\Delta)$. Then Proposition \ref{P:real-1} implies $\dim(M)=n(\Lg_\Delta)$.
We keep the notation in the proof of Proposition \ref{P:real-1}. Note that $\La_\Delta=\ker(e_1)$, $\La_\Delta^\bot=\bigcap_{i=2}^n\ker(e_i)$.
In this case, there exists $\lambda_1\in\cL_F(A_\Delta^\bot,\mu)$ such that $\dim E_F^{\lambda_1} = \dim E_F = n(\Lg_\Delta)$. In particular, we have $\cL_F(A_\Delta^\bot,\mu)=\{\lambda_1\}$.
Since
$$\{-2e_1,-e_1\pm e_2,\ldots,-e_1\pm e_n\}\subset\Sigma\sm\Sigma_Q\subset\RR_+\cL_F(A,\mu)$$
and
$$\{\lambda|_{\La_\Delta^\bot}:\lambda\in\cL_F(A,\mu)\}=\cL_F(A_\Delta^\bot,\mu)=\{\lambda_1\},$$
there exists $c>0$ such that
$$\cL_F(A,\mu)\cap\RR e_1=\{-ce_1\}$$ and
$$\cL_F(A,\mu)\cap\RR(-e_1\pm e_i)=\{-ce_1\pm ce_i\}, \qquad 2\le i\le n.$$
Let $\cL_F(A_\Delta,\mu)\subset\La_\Delta^*$ be the set of Lyapunov functionals for the $A_\Delta$-action on $E_F$ with respect to $\mu$. Then
$$\cL_F(A_\Delta,\mu)=\{\lambda|_{\La_\Delta}:\lambda\in\cL_F(A,\mu)\}\supset\{0,\pm ce_2,\ldots,\pm ce_n\}.$$
By Zimmer's cocycle superrigidity theorem, there is a representation $\phi:\Lg_\Delta\to\Lgl(n(\Lg_\Delta),\RR)$ whose set of restricted weights is equal to $\cL_F(A_\Delta,\mu)$ (hence contains $\{0,\pm ce_2,\ldots,\pm ce_n\}$), and such that for every $\eta\in\cL_F(A_\Delta,\mu)$, the dimension of the restricted weight space of $\phi$ with restricted weight $\eta$ is equal to
$$\dim\bigoplus_{\lambda\in\cL_F(A,\mu),\lambda|_{\La_\Delta}=\eta}E_F^\lambda.$$
Note that for the three cases of $G$ in Proposition \ref{P:real-2}, the subalgebra $\Lg_\Delta$ is isomorphic to $\Lsu(m-1,n-1)$, $\Lsp(m-1,n-1)$ or $\Lso^*(2m+2n-4)$, respectively. Since $\phi$ is nontrivial, by Lemma \ref{L:rep}, $\phi$ is equivalent to the standard representation of $\Lg_\Delta$. Note also that if we view $\Lg_\Delta$ as a subalgebra of $\Lgl(m+n-2,\CC)$ (resp. $\Lgl(m+n-2,\HH)$, $\Lgl(2m+2n-4,\CC)$) in the natural way and let $E_{i,j}$ denote the matrix of appropriate size with $1$ in the $(i,j)$-entry and $0$ elsewhere, then $\La_\Delta$ can be chosen such that its basis dual to $\{e_2,\ldots,e_n\}$ is
$$\{E_{m-n+i,m+n-i-1}+E_{m+n-i-1,m-n+i}:1\le i\le n-1\}$$
if $\Lg_\Delta=\Lsu(m-1,n-1)$ or $\Lsp(m-1,n-1)$ (see \cite[Section VI.4, Example 2]{Kn}), and is
$$\{E_{2i-1,m+n+2i-2}-E_{2i,m+n+2i-3}-E_{m+n+2i-3,2i}+E_{m+n+2i-2,2i-1}:1\le i\le n-1\}$$
if $\Lg_\Delta=\Lso^*(2m+2n-4)$ (see \cite[Problem VI.6]{Kn}).
We divide the remainder of the proof into two cases.

\textbf{Case 1.} Suppose $G$ is $\SU(n,n) \ (n\ge4)$, $\Sp(n,n) \ (n\ge3)$ or $\SO^*(4n) \ (n\ge4)$.
Then it follows from the above form of $\La_\Delta$ that $0$ is not a restricted weight of the standard representation of $\Lg_\Delta$, a contradiction.

\textbf{Case 2.} Suppose $G$ is $\SU(m,n) \ (m>n\ge3)$, $\Sp(m,n) \ (m>n\ge3)$ or $\SO^*(4n+2) \ (n\ge3)$.
Then the set of restricted weights of the standard representation of $\Lg_\Delta$ is $\{0,\pm e_2,\ldots,\pm e_n\}$. Moreover, if we let $V_0$ denote the restricted weight space of weight $0$, then $\dim V_0=\epsilon_G(m-n)$.
This implies that $c=1$. Thus
\ary \label{eq determinerLFAmucape1}
\{\lambda\in\cL_F(A,\mu):\lambda|_{\La_\Delta}=0\}=\cL_F(A,\mu)\cap\RR e_1=\{-e_1\}.
\eary
It follows that $\dim E_F^{-e_1}=\dim V_0=\epsilon_G(m-n)$.

Fix an arbitrary $a \in A$ wth $e_1(a) < 0$. By \eqref{eq determinerLFAmucape1}, we have
\begin{align*}
 (-e_1)(a) \dim E_{F}^{-e_1} & = (-e_1)(a) \epsilon_G(m-n) = (-e_1)(a) \dim \Lg_{-e_1}\\
 & <  (-e_1)(a) \dim \liealgg_{-e_1}  +  (-2e_1)(a) \dim \liealgg_{-2e_1}.
\end{align*}
By Corollary \ref{cor entropyinequalitytononatomic}, we deduce that $\mu^{G_{[-e_1]}}_x$ is non-atomic for $\mu$-a.e. $x$.

For $\mu$-a.e. $x\in M^\alpha$, the non-atomicity of $\mu^{G_{[-e_1]}}_x$ implies that there exists a nonzero vector $\bx_1+\bx_2\in\Lg_{[-e_1]}=\Lg_{-e_1}\oplus\Lg_{-2e_1}$, where $\bx_1\in\Lg_{-e_1}$ and $\bx_2\in\Lg_{-2e_1}$, such that $\exp(\bx_1+\bx_2)\in\supp \mu^{G_{[-e_1]}}_x$. Let us define $\beta,\gamma\in\Sigma$ and $\bx\in\Lg_\beta\sm\{0\}$ as follows:
\begin{itemize}
  \item If $\bx_1\ne0$, we let $\beta=-e_1$, $\gamma=e_2$ and $\bx=\bx_1$;
  \item If $\bx_1=0$, then $\bx_2\ne0$, and we let $\beta=-2e_1$, $\gamma=e_1+e_2$, $\bx=\bx_2$.
\end{itemize}
In each case, we have $\gamma\in\Sigma_Q$, $\beta+\gamma=-e_1+e_2\notin\Sigma_Q$, and $\langle\beta,\beta+2\gamma\ra\ge0$. So by Lemma \ref{L:root-space}, we have $[\bx,\Lg_\gamma]=\Lg_{\beta+\gamma}$. It follows from $\gamma\in\Sigma_Q$ that $\mu^{G_{[\gamma]}}_x$ is a Haar measure on $G_{[\gamma]}$. We then  deduce from Proposition \ref{lem noncommutinginvariance} that $\mu^{G_{\beta+\gamma}}_x$ is a Haar measure on $G_{\beta+\gamma}$. Thus $\mu$ is invariant under $G_{\beta+\gamma}$, which contradicts the fact that $\beta+\gamma\notin\Sigma_Q$.
\end{proof}

\begin{proof}[Proof of \eqref{E:s-uni1}]
Let $G=\SU(n,n)$ or $\SO^*(4n)$ with $n\ge2$. We need to prove that if $\alpha$ is an action of a lattice $\Gamma$ in $G$ on a compact manifold $M$ such that there is a $P$-invariant $P$-ergodic Borel probability measure $\mu$ on $M^\alpha$ that is not $G$-invariant, then
$$\dim(M)\ge\min\{2\epsilon_G(n-1)+1,n^2\}.$$
Let $Q=\Stab_G(\mu)$. Then $\Pi_Q\ne\Pi$. Let $k\in\{1,\ldots,n\}$ be the largest index such that $\alpha_k\notin\Pi_Q$. We divide the argument into 3 cases.

\textbf{Case 1.} Suppose $n\ge4$ and $k=1$. Let $\Delta=\Pi_Q=\{\alpha_2,\ldots,\alpha_n\}$. If $G$ is $\SU(n,n)$ (resp. $\SO^*(4n)$), then $\Lg_\Delta$ is isomorphic to $\Lsu(n-1,n-1)$ (resp. $\Lso^*(4n-4)$). In each case we have $n(\Lg_\Delta)=2\epsilon_G(n-1)$. By Proposition \ref{P:real-2}, we have
$$\dim(M)\ge n(\Lg_\Delta)+1=2\epsilon_G(n-1)+1.$$

\textbf{Case 2.} Suppose $n\ge5$ and $2\le k\le n-3$. Let $\Delta=\{\alpha_{k+1},\ldots,\alpha_n\}$. If $G$ is $\SU(n,n)$ (resp. $\SO^*(4n)$), then $\Lg_\Delta$ is isomorphic to $\Lsu(n-k,n-k)$ (resp. $\Lso^*(4n-4k)$). By Proposition \ref{P:real-1}, we have
$$\dim(M)\ge kn(\Lg_\Delta)=2\epsilon_Gk(n-k)\ge 4\epsilon_G(n-2)\ge 2\epsilon_G(n-1)+1.$$

\textbf{Case 3.} Suppose the assumptions in Cases 1 and 2 do not hold. Then either $n\in\{2,3\}$, or $n\ge4$ and $k\ge n-2$. In particular, we always have $k\ge n/3$. It follows from $\alpha_k\notin\Pi_Q$ that $\Sigma_Q$ does not contain the following $k(4n-3k+1)/2$ roots:
$$\{-2e_i:1\le i\le k\}\cup\{-e_i-e_j:1\le i<j\le k\}\cup\{-e_i\pm e_j:1\le i\le k, k+1\le j\le n\}.$$
By Proposition \ref{prop dim geq r0} and \eqref{E:r0}, we have
$$\dim(M)\ge r_0(Q)\ge k(4n-3k+1)-k\ge n^2.$$
This proves \eqref{E:s-uni1}.
\end{proof}

\begin{proof}[Proof of \eqref{E:s-uni2}]
Let $G$ be one of the groups in \eqref{E:s-uni2}. We need to prove that if $\alpha$ is an action of a lattice $\Gamma$ in $G$ on a compact manifold $M$ such that there is a $P$-invariant $P$-ergodic Borel probability measure $\mu$ on $M^\alpha$ that is not $G$-invariant, then
$$\dim(M)\ge\min\{\epsilon_G(m+n-2)+1,n(n+1)\}.$$
Let $Q=\Stab_G(\mu)$. Then $\Pi_Q\ne\Pi$. Let $k\in\{1,\ldots,n\}$ be the largest index such that $\alpha_k\notin\Pi_Q$. We divide the argument into 3 cases.

\textbf{Case 1.} Suppose $n\ge3$ and $k=1$. Let $\Delta=\Pi_Q=\{\alpha_2,\ldots,\alpha_n\}$. If $G$ is $\SU(m,n)$ (resp. $\Sp(m,n)$, $\SO^*(2m+2n)$), then $\Lg_\Delta$ is isomorphic to $\Lsu(m-1,n-1)$ (resp. $\Lsp(m-1,n-1)$, $\Lso^*(2m+2n-4)$). In each case we have $n(\Lg_\Delta)=\epsilon_G(m+n-2)$. By Proposition \ref{P:real-2}, we have
$$\dim(M)\ge n(\Lg_\Delta)+1=\epsilon_G(m+n-2)+1.$$

\textbf{Case 2.} Suppose $n\ge4$ and $2\le k\le n-2$. Let $\Delta=\{\alpha_{k+1},\ldots,\alpha_n\}$. If $G$ is $\SU(m,n)$ (resp. $\Sp(m,n)$, $\SO^*(2m+2n)$), then $\Lg_\Delta$ is isomorphic to $\Lsu(m-k,n-k)$ (resp. $\Lsp(m-k,n-k)$, $\Lso^*(2m+2n-4k)$). By Proposition \ref{P:real-1}, we have
$$\dim(M)\ge kn(\Lg_\Delta)=\epsilon_Gk(m+n-2k)\ge 2\epsilon_G(m+n-4)\ge\epsilon_G(m+n-2)+1.$$

\textbf{Case 3.} Suppose the assumptions in Cases 1 and 2 do not hold. Then $k\ge n-1$. It follows from $\alpha_k\notin\Pi_Q$ that $\Sigma_Q$ does not contain the following $k(4n-3k+1)/2$ roots:
$$\{-2e_i:1\le i\le k\}\cup\{-e_i-e_j:1\le i<j\le k\}\cup\{-e_i\pm e_j:1\le i\le k, k+1\le j\le n\}.$$
By Proposition \ref{prop dim geq r0} and \eqref{E:r0}, we have
$$\dim(M)\ge r_0(Q)\ge k(4n-3k+1)\ge n(n+1).$$
This proves \eqref{E:s-uni2}.
\end{proof}

\appendix

\section{The quantities $v(G)$, $v(G_{\mathrm{cpt}})$, $n(G)$ and $r(G)$}\label{S:vnd}

{\footnotesize
\begin{table}[h]\hspace{-2cm}
\begin{tabular}{|c|c|c|c|c|}\hline
$G$ & $v(G)$ & $v(G_{\mathrm{cpt}})$ & $n(G)$ & $r(G)$ \\  \hline
$\SL(n,\CC), \ n\ge3$ & $2n-2$ & $\begin{matrix} 2n-2, & n\ne4; \\ 5, & n=4 \end{matrix}$ & $2n$ & $n-1$ \\  \hline
$\Sp(2n,\CC), \ n\ge3$ & $4n-2$ & $4n-4$ & $4n$  & $2n-1$ \\  \hline
$\SO(n,\CC), \ n\ge7$ & $2n-4$ & $n-1$ & $2n$  & $n-2$ \\  \hline
$\SL(n,\HH), \ n\ge3$ & $4n-4$ & $4n-2$ & $4n$  & $n-1$ \\  \hline
$\SO^+(m,n), \ m> n\ge2$ & $m+n-2$ & $m+n-1$ & $m+n$  & $2n-1$ \\  \hline
$\SU(m,n), \ m\ge n\ge2$ & $\begin{matrix} 2m+2n-3, & (m,n)\ne(2,2);  \\ 4, & (m,n)=(2,2) \end{matrix}$ & $\begin{matrix} 2m+2n-2, & (m,n)\ne(2,2); \\ 5, & (m,n)=(2,2) \end{matrix}$ &
$\begin{matrix} 2m+2n, & (m,n)\ne(2,2); \\ 6, & (m,n)=(2,2) \end{matrix}$  & $2n-1$ \\  \hline
$\Sp(m,n), \ m\ge n\ge2$ & $\begin{matrix} 4m+4n-5, & (m,n)\ne(2,2);  \\ 10, & (m,n)=(2,2) \end{matrix}$ & $4m+4n-4$ & $4m+4n$  &  $2n-1$ \\  \hline
$\SO^*(2n), \ n\ge4$ & $\begin{matrix} 4n-7, & n\ne 4,6;  \\ 6, & n=4; \\ 15, & n=6 \end{matrix}$ & $2n-1$ & $\begin{matrix} 4n, & n\ge5; \\ 8, & n=4 \end{matrix}$  & $2[\frac{n}{2}]-1$ \\  \hline
$E_6$ & $32$ & $26$ & $54$  & $16$ \\  \hline
$E_7$ & $54$ &  $54$ & $112$  & $27$ \\  \hline
$E_8$ & $114$ & $112$ & $496$  & $57$ \\  \hline
$F_4$ & $30$ &  $16$ & $52$  & $15$ \\ \hline
$G_2$ & $10$ &  $6$ & $14$  & $5$ \\ \hline
\end{tabular}
\caption{The quantities $v(G), v(G_{\mathrm{cpt}}), n(G)$ and $r(G)$ for some non-split simple Lie groups $G$. The notations $E_n, F_4, G_2$ refer to \emph{complex} simple Lie groups of the corresponding types.}
\label{Tb1}
\end{table}
}

\
\



\end{document}